\documentclass[12pt]{article}
\newsavebox{\toy}
\savebox{\toy}{\framebox[0.65em]{\rule{0cm}{1ex}}}
\newcommand{\QED}{\usebox{\toy}\end{demo}}
\newenvironment{property}%
{\begin{list}{}{\setlength{\rightmargin}{0pt}%
\setlength{\itemsep}{0pt}}}{\end{list}}
\newlength{\templength}
\newcommand{\bp}{\setlength{\templength}{\labelwidth}%
\setlength{\labelwidth}{2em}\begin{property}}
\newcommand{\ep}{\end{property}\setlength{\labelwidth}{\templength}}
\newtheorem{theorem}{Theorem}[subsection]
\newtheorem{lemma}[theorem]{Lemma}
\newtheorem{proposition}[theorem]{Proposition}
\newtheorem{corollary}[theorem]{Corollary}
\newtheorem{conjecture}[theorem]{Conjecture}
\newtheorem{definition}[theorem]{Definition}
\newtheorem{example}[theorem]{Example}
\newtheorem{remark}[theorem]{Remark}
\newtheorem{exercise}{Exercise}[subsection]
\newtheorem{assumption}{Assumption}
\newcommand{\Thm}[1]{Theorem \ref{Thm.#1}}
\newcommand{\Lem}[1]{Lemma \ref{Lem.#1}}

\newcommand{\Prop}[1]{Proposition \ref{Prop.#1}}
\newcommand{\Cor}[1]{Corollary \ref{Cor.#1}}

\newcommand{\Theorem}[1]{\begin{theorem}\label{Thm.#1}}
\newcommand{\Lemma}[1]{\begin{lemma}\label{Lem.#1}}
\newcommand{\Proposition}[1]{\begin{proposition}\label{Prop.#1}}
\newcommand{\Corollary}[1]{\begin{corollary}\label{Cor.#1}}
\newcommand{\Assumption}[1]{\begin{assumption}\label{Ass.#1}\rm}
\newcommand{\Definition}[1]{\begin{definition}\label{Def.#1}\rm}
\newcommand{\Remark}[1]{\begin{remark}\label{Rem.#1}\rm }
\newcommand{\Exercise}[1]{\begin{exercise}\label{Exe.#1}\rm }
\newcommand{\Example}[1]{\begin{example}\label{Exa.#1}\rm }
\newcommand{\bd}{\begin{displaymath}}
\newcommand{\ed}{\end{displaymath}}
\newcommand{\bdn}{\begin{equation}}
\newcommand{\bdnl}{\begin{equation}\label}
\newcommand{\edn}{\end{equation}}
\newcommand{\barray}{\begin{array}}
\newcommand{\earray}{\end{array}}
\newcommand{\bds}{\begin{description}}
\newcommand{\eds}{\end{description}}
\newcommand{\bitemize}{\begin{itemize}}
\newcommand{\eitemize}{\end{itemize}}
\newcommand{\benumerate}{\begin{enumerate}}
\newcommand{\eenumerate}{\end{enumerate}}
\newcommand{\btabbing}{\begin{tabbing}}
\newcommand{\etabbing}{\end{tabbing}}
\newcommand{\bcenter}{\begin{center}}
\newcommand{\ecenter}{\end{center}}
\newcommand{\bflushright}{\begin{flushright}}
\newcommand{\bflushleft}{\begin{flushleft}}
\newcommand{\eflushright}{\end{flushright}}
\newcommand{\eflushleft}{\end{flushleft}}
\newcommand{\bdnn }{\begin{eqnarray*}}
\newcommand{\ednn }{\end{eqnarray*}}
\newcommand{\bdmn}{\begin{eqnarray}}
\newcommand{\edmn}{\end{eqnarray}}

\newcommand{\nn}{\nonumber}
\newcommand{\SSC}[1]{\section{#1}\setcounter{equation}{0}}

\newcounter{biblio}
\newenvironment{references}%
{\begin{list}{[\arabic{biblio}]}{\usecounter{biblio}%
\setlength{\leftmargin}{2.5em}\setlength{\rightmargin}{0pt}%
\setlength{\labelwidth}{2em}\setlength{\itemsep}{0pt}}}{\end{list}}
\newcommand{\References}%
{\vspace{2.8ex plus .3ex minus .3ex}%
\begin{center}{\bf References}\end{center}\begin{references}}

\usepackage{latexsym}
\usepackage{graphicx}
\usepackage{amssymb}
\usepackage{mathrsfs}
\newcommand{\R}{{\mathbb{R}}}
\newcommand{\rd}{\R^d}

\newcommand{\ra }{\rightarrow }
\newcommand{\lra }{\longrightarrow }

\newcommand{\ov}{\overline}
\newcommand{\tl}{\widetilde}

\newcommand{\Llra}{\Longleftrightarrow }

\newcommand{\vvs}{\vspace{2ex}}
\newcommand{\vs}{\vspace{1ex}}

\newcommand{\lef}{\left}
\newcommand{\rig}{\right}
\newcommand{\ri}{\right}
\newcommand{\st}{\stackrel}

\newcommand{\8}{\infty}

\newcommand{\6}{\partial}
\newcommand{\dps}{\displaystyle}
\newcommand{\sub}{\subset}

\newcommand{\bsh}{\backslash}

\newcommand{\half}{\mbox{$\frac{1}{2}$}}

\newcommand{\inflim}{\mathop{\underline{\lim}}}
\newcommand{\suplim}{\mathop{\overline{\lim}}}
\newcommand{\epty}{\emptyset}
\renewcommand{\a}{\alpha}
\renewcommand{\b}{\beta}
\newcommand{\gm}{\gamma}
\newcommand{\Gm}{\Gamma}

\newcommand{\del}{\delta}

\newcommand{\e}{\varepsilon}

\newcommand{\z}{\zeta}
\newcommand{\h}{\eta}

\newcommand{\lm}{\lambda}

\newcommand{\m}{\mu}
\newcommand{\n}{\nu}

\newcommand{\rh}{\rho}

\newcommand{\s}{\sigma}

\newcommand{\om}{\omega}
\newcommand{\W}{\Omega}

\newcommand{\cB }{{\cal B}}
\newcommand{\cC }{{\cal C}}
\newcommand{\cD }{{\cal D}}

\newcommand{\cG }{{\cal G}}

\newcommand{\cK }{{\cal K}}
\newcommand{\cL }{{\cal L}}
\newcommand{\cM }{{\cal M}}

\newcommand{\cO }{{\cal O}}

\newcommand{\zk}{\zeta^{(k)}}
\newcommand{\pk}{\psi^{(k)}}
\newcommand{\vak}{\varphi^{(k)}}

\newcommand{\Yt}{{\mathsf Y}^{(t)}}

\newcommand{\Crit}{\tt Crit}

\makeatletter
\def\section{\@startsection{section}{1}{\z@}{-3.5ex plus -1ex minus 
 -.2ex}{2.3ex plus .2ex}{\bf}}
\def\subsection{\@startsection{subsection}{2}{\z@}{-3.25ex plus -1ex minus 
 -.2ex}{1.5ex plus .2ex}{\bf}}
\makeatother

\begin{document}

\bcenter

\large{\bf 
Localization Transition for Polymers 
in Poissonian Medium}\footnote{December 6, 2012}

\vvs

\vvs \normalsize
\noindent Francis COMETS
\footnote{
Partially
supported by  CNRS (UMR 7599
Probabilit{\'e}s et Mod{\`e}les
Al{\'e}atoires)
}
\\

\vs \small
Universit{\'e} Paris Diderot - Paris 7, \\
Math{\'e}matiques, Case 7012\\
75205 Paris cedex 13, France \\
email: {\tt comets@math.univ-paris-diderot.fr} \\
{\tt http://www.proba.jussieu.fr/$\sim$comets}

\vvs \normalsize

\noindent Nobuo YOSHIDA\footnote{Partially
supported by JSPS Grant-in-Aid for Scientific Research, 
Kiban (C) 17540112}\\

\vs \small
Division of Mathematics \\
Graduate School of Science \\
Kyoto University,\\
Kyoto 606-8502, Japan.\\
email: {\tt nobuo@math.kyoto-u.ac.jp}\\
{\tt http://www.math.kyoto-u.ac.jp/$\sim$nobuo/}

\ecenter

\vs
\begin{abstract}
We study a model of directed polymers in random environment in dimension $1+d$, given
by a Brownian motion in a Poissonian potential. We study the effect of the density and the 
strength of inhomogeneities, respectively the intensity parameter $\n$ of the Poisson field
and  the temperature inverse $\b$. Our results are: (i) fine information on the phase diagram, with
quantitative estimates on the critical curve; (ii) pathwise localization at low temperature and/or
large density; (iii) complete localization in a favourite corridor for large 
$\n \b^2$ and bounded $\b$.
\end{abstract}
\footnotesize

\footnotesize
\noindent{\bf Short Title:} Brownian Polymers



\tableofcontents

\vspace{1cm}

\normalsize

\SSC{Introduction}
We study in the present article the long-time behavior of the Brownian directed polymer in
dimension $d \geq 1$, under the influence of a random Poissonian environment, 
as introduced in \cite{CY05}. Let $B=(B_t)_{t \geq 0}$ be the canonical Brownian motion in
$\R^d$ starting at the origin, and $P$ the Wiener measure. Let also $Q$ be the law of  the
canonical 
Poisson point process $\eta$ in $\R_+ \times \R^d$ with intensity measure $\nu dt dx$, where $\n$ is a positive parameter. 
Denote by $U(x) \sub \rd$ the closed ball with the unit 
volume, centered at $x \in \rd$, and by 
%
$$
H_t^\h(B) = {\rm Card}\lef\{ {\rm points}\;(s,x)\; {\rm in \ } \h : 
s \leq t, \; x \in U(B_s)\ri\},
$$
the number of Poisson points which are "seen" by the path $B$ up to time $t$. We are interested
in the behavior of $B$ for large $t$ and typical $\h$, under the polymer  measure $\mu_t$
on the path space $C(\R_+ \ra \rd)$ given by
\begin{equation}
d\mu_t = \frac{1}{Z_t}  \exp \{ \b H_t^\h (B) \} \; dP
\end{equation}
Here, $Z_t$ is the normalizing constant, and $\b \in \R$ is a parameter. Its absolute value, being
proportional to the temperature inverse,  measures the strength of the inhomogeneities 
produced by the random environment, whereas its sign indicates whether the path prefer
to visit Poisson points  or to stay clear from them. The model has an interpretation in terms of branching
Brownian motions in random environment. When $\b>0$, the  Poisson points  are catalyzers,
and every point $(s,x)$ in $\eta$ causes an instantaneous (possibly multiple) branching with mean  offspring $e^\b$
to every individual passing by $x$ within distance $r_d$ at time $s$. 
When $\b<0$, the  Poisson points  are soft obstacles,
and every point $(s,x)$ in $\eta$ kills individual passing by $x$ within distance $r_d$ 
at time $s$ with probability $1-e^\b$. Then, $Z_t$ is the average population at time $t$ in the environment $\eta$,  e.g., the survival probability when $\b \leq 0$, and the restriction of $\mu_t$ 
to time interval $[0,t]$ is the law of the ancestral line of a randomly selected individual 
in the population at time $t$, conditionally on survival. 

This particular model was considered in \cite{CY05, CY04} with $\n=1$.  As is the rule for general polymer models, a salient feature in dimension
$d \geq 3$, is a phase transition between a high temperature phase ($|\b|$ small) where 
inhomogeneities are 
inessential, and a low temperature phase ($|\b|$ large) where 
inhomogeneities are crucial. The phases are first defined by thermodynamic functions, namely by the dis/agreement
of the quenched and annealed free energies, 
but it was discovered that they correspond to delocalized and localized behavior respectively,
see \cite{CaHu02, CSY03}  for simpler models. The thermodynamic transition coincides with localization transition. 
The counterpart of the diffusive and of localized  behaviors for branching Brownian motions  have been studied in \cite{Shio1, Shio2} with a particular dynamics.

The proof that there exists a high temperature phase goes back to \cite{ImSp88, Bol89}.
The sub-region defined by the condition in (\ref{eq:L2}) below, is called the $L^2$-region, 
is where $Z_t/Q Z_t$ is bounded in $L^2$. There, second moment method works, showing that the polymer measure is 
much similar to the Wiener measure.
Non perturbative results covering the full high temperature region are rare.

At  low temperature, localization properties are traditionally formulated for the end point of the polymer  as in \cite{CaHu02, CSY03}.
Recent results in \cite{CC11} for the parabolic Anderson model, have led to substantial progress in understanding that localization holds in a stronger, pathwise manner there:
The polymer path
spends a positive proportion of time at the same location as some particular path depending on the realization of the medium.  The proof 
crucially uses the Gaussian nature of the environment, via integration by parts formula,
and it is not clear how general pathwise localization  is.
We observe that the concentration effect is a global phenomenon in our model, in contrast with heavy-tails potentials where only extreme statistics are relevant \cite{AuLo10, HaMa07} and it
only matters that the random path visits the few corresponding locations.
By nature, available information in the present model concerns replica overlaps. 
Moments,  covariances, conditional moments can often be represented as expected values of independent copies 
of the paths sharing the same environment, so-called  replica.
A necessary step is to extract, from the latter,
information on a single polymer.

The present model is quite natural, and interestingly enough, it is related to other polymer models.
For instance, the mean field limit, $\nu \to \8$ and $\beta \to 0$ in such a way that $\nu \beta^2 \to b^2 >0$,
is the Brownian directed polymer in a Gaussian environment. In this model introduced in
\cite{RoTi05},  the environment is the generalized Gaussian process $g(t,x)$ with mean 0 and
covariance
$$
Q[g(t,x)g(s,y)] 
= b^2 \delta (t-s) |U(x) \cap U(y)|,
$$
where $|\;\cdot \; |$ above denotes the Lebesgue measure. 
As mentioned in \cite{La11}, the proofs of superdiffusivity in one space dimension and
the analysis of the influence of spatial correlations for the Gaussian environment case, can be
adapted to our case of a Poissonian environment.

A few exactly solvable polymer models 
are known so far, all of them being for $d=1$: (i) the infinite series of Brownian queues \cite{OCYor01, MoOc07}, which is 
a limit of strongly asymmetric polymers \cite{Moreno10}; (ii) the discrete model with log-gamma weights
\cite{Sep09, GSep12};  
(iii) the Hopf-Cole solution of the one-dimensional Kardar-Parisi-Zhang
equation \cite{AmCoQu11}, which is expected the universal scaling limit for polymers \cite{Corwin11}.
Exact solutions are also available at zero temperature, via determinantal processes.   
Note that, from \cite{Lac10}, there is no high temperature region in dimension $d=1$ and $2$.

For disordered polymer pinning on an interface, which also shows  localized and a delocalized phases, 
estimating the critical curve is an important and difficult problem \cite{GB, denHollander-SF}. 
The similar remark also holds for bulk disorder with long range correlation \cite{CaTiVi08}.
Unrelated disordered systems,  including the Sherrington-Kirkpatrick model, have seen the emergence of smart 
interpolation techniques \cite{Tala, Guerra03, GuerraToni}, allowing to compare the free energy
with that of an auxiliary, simpler model.

We finally mention the relations to the non directed, model of Brownian motion in a space-dependent (but time independent) Poissonian potential, which has been extensively studied 
in many different perspectives; We refer to
\cite{Szn98} for a detailed overview. 
The spectral theory approach and the coarse-graining method of enlargement of obstacles 
developed in \cite{Szn98} does not apply to our directed model. The latter one can be 
thought as the case of very strong drift  in a fixed direction, or, equivalently, the problem of long crossings, which has been considered in such continuous models \cite{Wut98b} as well as in discrete
ones \cite{IoVe12, Zyg12}.
\medskip

The main objective of the present paper is to study the joint effect of the density and the strength of 
inhomogeneities, i.e., the influence of the intensity parameter $\n$ and of the temperature inverse 
$\b$. We obtain qualitative and quantitative estimates on the critical curve separing the two phases in the plane $(\b,\n)$. We find some auxiliary curves in this plane  along which the difference between the annealed and quenched free energies
is monotone, thus they do not re-enter a phase after leaving. In the spirit of the interpolation techniques mentioned above, the control of the sign of the derivative is made possible in the present model by the integration by parts formula for the Poisson process. 
A second set of results is for the pathwise localization in the localized phase, it applies in all space dimension $d =1,2\ldots$ (the polymer "physical" dimension being $1+d$).
We define a trajectory -- we call it the favourite path for obvious reasons -- depending on the realization of the medium and on the parameters, in the vicinity of  which the random polymer path
spends a positive fraction of time. Moreover, the localization in the favourite corridor 
becomes complete in some region of the parameter space: 
we show that this fraction converges to
1 as $\n \b^2 \to \8$ with $\b$ remaining bounded. To parallel the notion of geodesics
in last passage percolation \cite{New95}, the favourite path can be
viewed as a "fuzzy geodesics", and this one is essentially unique in this asymptotics.
\medskip

Our paper is organized as follows: We start with notations and previous results. We then formulate
our main results in section 3. In the next section we define the favourite path, the overlap between two polymer paths (replica), and prove a "two-to-one lemma", which extracts information on a single polymer path from the overlap of two replica. Section 5 contains the proofs of localization,
except for an estimate, needed for complete localization, of the discrepancy of quenched and annealed free energies, which is obtained in section 6.
The final section is devoted to the estimates of the critical curve.

\SSC{The model of Brownian directed polymers in random environment}
\label{themodel}
\subsection{Preliminaries}

We set some more notations. The environment
$\h $ is the Poisson random measure on $\R_+ \times \rd$ 
with the intensity $\n>0$, 
defined on the probability space $(\cM, \cG, Q)$, 
with $\cM$ is the set of integer-valued 
Radon measure on $\R_+ \times \rd$, 
$\cG$ is the $\sigma$-field  generated by the variables $\h (A)\; , \; A \in \cB (\R_+ \times \rd)$.
$Q$ is the unique probability measure on $(\cM, \cG)$ such that,
for disjoint and bounded $A_1, ..., A_n \in \cB (\R_+ \times \rd)$,
the variables  $\h (A_j)$ are independent with Poisson distribution of mean
$\n |A_j|$; 
Here, 
$ | \cdot |$ denotes the Lebesgue measure 
on $\R^{1+d}$.
For $t>0$, it is natural and convenient to introduce its restriction
\bdnl{ht}
\h_t(A)  =  \h (A \cap ((0,t] \!\times \! \rd))\;,\; \; \; 
A \in \cB (\R_+ \times \rd).
\edn
We denote by  $V_t$  the 
tube around the 
graph
$\{(s,B_s)\}_{0 < s \leq t}$ of the Brownian path,
\bdnl{Vt}
V_t=V_t(B)=\{ (s,x)\; ; \; s \in (0,t], \; x \in U (B_s)\},
\edn
where $U(x) \sub \rd$ is the closed ball with the unit 
volume, centered at $x \in \rd$. ($U(x)$ has radius $r_d$.)
Then, for any $t>0$,
the polymer  measure $\m_t^x$ can be expressed  as
\bdnl{mnen}
d \m_t ={(Z_t)^{-1}} {\exp \lef( \b \h (V_t)\rig)}\; dP,
\edn
with the partition function $Z_t$ 
\bdn 
Z_t=P[\exp \lef( \b \h (V_t)\rig)]\;. \label{Zt}
\edn

Let $f,g:I \ra (0,\8)$ be functions and $a \in \ov{I}$, 
where $I \sub \R$ is an interval. We write 
$f(x) \sim g(x)$ ($x \ra a$), if $\lim_{x \ra a}f(x)/g(x)=1$. 
We write 
$f(x) \asymp g(x)$ ($x \ra a$), if 
$0<\inflim_{x \ra a}f(x)/g(x) \le \suplim_{x \ra a}f(x)/g(x) <\8$.

\subsection{Former results}
Denote by $\lambda$
the logarithmic moment generating function of a 
mean-one Poisson distribution,
\bdnl{lambda}
\lm = \lm (\b)=e^\b-1 \in (-1,\8)\;.
\edn
The quenched free energy $p_t(\b, \nu)$ of the polymer model with finite time horizon $t$  is
\bdnl{pt(b,nu)}
p_t(\b, \nu)=\frac{1}{t}Q \ln Z_t,
\edn
though $\frac{1}{t} \ln Q Z_t=\n \lm(\b)$ is  the annealed  free energy.
The case  of a fixed $\nu=1$ was considered in the papers \cite{CY05, CY04}, but the results trivially 
extend to a general $\nu$. We summarize them without repeating the proof.
\Theorem{th:psi} 
Let $d \ge 1, \nu >0$ and $\b \in \R$ be arbitrary. 
\bds
\item[(a)]
There exists a deterministic number $p(\beta, \nu) \in \R$ such that 
\bdmn
p(\b , \nu)&=&\lim_{t \nearrow \8}p_t(\b , \nu), \label{psias} \\
&=&\lim_{t \nearrow \8}\frac{1}{t}\ln Z_t,\; \; \; 
\mbox{$Q$-a.s. and in $L^2(Q)$.} \label{psias2}
\edmn
\item[(b)] 
The function $\b \mapsto p(\beta, \nu)$ is convex on $\R$, with
\bdnl{psile}
\nu \beta \leq  p(\b , \nu) \leq \nu \lambda .
\edn
\noindent
The function $\beta \mapsto\nu \lm (\beta) - p(\beta, \nu)$ 
is non-decreasing on $\R_+$ and  non-increasing on $\R_-$.
\item[(c)] 
There exist critical values $\b_c^\pm=\b_c^\pm(\n )=\b_c^\pm(d,\n )$ with 
$-\8 \leq \b^-_c \leq 0 \leq \b^+_c < +\8,$
such that 
\bdmn
p(\beta, \nu)=\nu \lm
\; \; & \mbox{if} &\; \; \b \in [\b_c^-, \b^+_c]\cap \R, 
\label{[b-,b+]} \\
p(\beta, \nu)< \nu \lm,
\; \; & \mbox{if} &\; \; \b \in \R \bsh [\b_c^-, \b^+_c].
\label{Rbsh[b-,b+]} 
\edmn
\item[(d)] 
For $d \ge 3$, $\b^-_c(d,\n)<0 <\b^+_c(d,\n)$, 
${\dps \lim_{d \nearrow \8}\b_c^\pm(d,\n)=\pm \8}$,  
and there exists $\n_c \in [1,\8]$ (cf. \Prop{prop:nu_c} below) such that
\bdnl{nuc}
\b^-_c(d,\n)\lef\{ \barray{ll}
=-\8 & \mbox{if $\n <\n_c$}, \\
>-\8 & \mbox{if $\n >\n_c$}.
\earray \ri.
\edn
More precisely, letting
$$a_{L^2} = \sup \lef\{ a>0:
 P\Big[ \exp \lef( {a \over 2}
\int^\8_0| U(0) \cap U(B_s ) |ds \ri)\Big] < \8 \ri\}>0,$$ then
 \begin{equation} \label{eq:L2}
\n \lm (\b)^2 < a_{L^2} \Longrightarrow p(\b, \n)=\n \lm (\b)\;,
\end{equation}
 and thus, $\n_c \geq a_{L^2}$ and
\bdnl{bc*sand*rough}
\b^-_c (d,\n) \le \ln \lef( 1-\sqrt{a_{L^2}/\n}\ri) 
< \ln \lef( 1+\sqrt{a_{L^2}/\n}\ri) \le \b^+_c (d,\n).
\edn
%
%
%
%
\eds 
\end{theorem}
For completeness, we mention a numerical lower bound for  $a_{L^2}$, and thus for $\nu_c$ itself, 
which can be derived from the techniques of section 4.2 in \cite{CY05}.
Let $\gm_d$ denote the smallest positive zero of the Bessel function 
$
J_{\frac{d-4}{2}} (\gm )  =  (\gm /2)^{{\frac{d-4}{2}} }\sum_{k \ge 0}
{(-\gm^2 /4)^k}/({k!\Gm ({\frac{d-4}{2}} +k+1)}),$ $\gm \ge 0.
$
Then, with $r_d$ the radius of the ball $U(0)$ with unit volume,
\begin{equation} 
\nu_c \geq \left( \frac{\gamma_d}{2r_d}\right)^2.
\end{equation}
The lower bound has value 1.265\ldots for $d=3$, 1.792\ldots for $d=4$, 2.190\ldots for $d=5$, and
$\liminf_{d \to \8} d^{-1/2} \n_c \geq \sqrt{e/(8\pi)}=0.329\ldots$.

\SSC{Main results} \label{sec:results}
\subsection{Phase diagram}

The parameter space $(0,\8)\times \R$ splits into two regions,
$$
\cD=\{(\n, \b): p(\nu,\b) = \n \lm(\b)\}, \qquad \cL=\{(\n, \b): p(\nu,\b) < \n \lm(\b)\}=\cD^c, 
$$
which are called {\it high temperature / low density region} and {\it low temperature / high density region} respectively. 
The name 
is justified by observing that infinite temperature, or equivalently, $\b =0 $,
 belongs to $\mathcal D$, though  zero density $\n=0$ belongs to this set. 
We already know from \cite{CY05} that they correspond to end-point delocalized and localized phase, see (\ref{eq:pathloc1}) below.
In the next section, we will discuss deeper aspects of localization.

\begin{figure} \centering
\includegraphics[width=15cm,height=100mm]{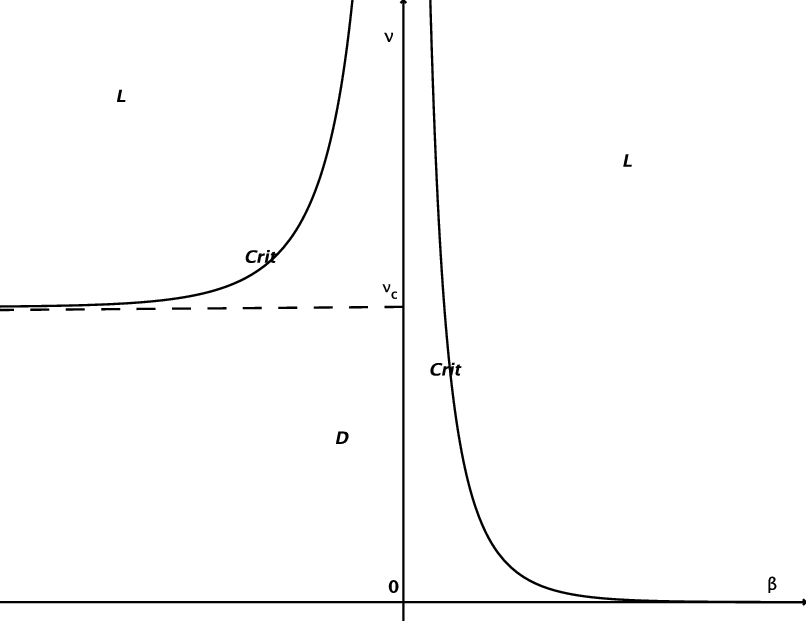}
\caption{Phase diagram, $d \geq 3$. The high temperature/low density and low temperature/high density phases $\cal D$ and $\cal L$ are separated by the critical curve $\{(\beta_c^+(\nu), \nu); \nu >0\} \cup \{(\beta_c^-(\nu), \nu);
\nu > \nu_c\}$} 
\label{f-diag}
\end{figure}

We state some properties of $\cL, \cD$, and of the critical curve separating the two sets,
$$
\Crit=\cD \cap \overline \cL .
$$
It is proved in \cite{bertin}
that $\cD$ reduces to the semi-axis $\{\nu>0, \b=0\}$ 
in dimensions $d=1$ and $d=2$, so we focus on the case of 
a larger dimension. By \Thm{th:psi} , (\ref{[b-,b+]}), (\ref{Rbsh[b-,b+]}),
we already know that $\Crit$ is the union of the 
graphs of the functions $\nu \mapsto
\b_c^-(d,\n)$ and $\nu \mapsto \b_c^+(d,\n)$. Moreover, by \Lem{mono*nu}
we see that $ \n \lm(\b)-p(\nu,\b)$ is non-decreasing in $\nu$,  non-decreasing in $\b$ for $\b \geq 0$ and non-increasing in $\b$ for $\b \leq 0$.
The qualitative features of the phase diagram are summarized in figure \ref{f-diag},
corresponding to statements all through the present Section \ref{sec:results}.

We first answer some questions which were left open in  \Thm{th:psi}, (e). 
\Proposition{prop:nu_c}
For all dimension $d \geq 3$, we have 
\bdnl{n_c<8}
 \nu_c <\8
\edn
and 
\bdnl{eq:betanuc}
\beta_c^-(\nu_c)=-\8,
\edn
\end{proposition} 
Recall for completeness that, in lower dimensions $d=1,2,$ we know from \cite{bertin} that $\beta^-_c(\nu)=0$ for all $\nu$.
We prove (\ref{n_c<8}) as a part of \Cor{decay*of*b_c} below. 
By (\ref{n_c<8}), the set $\cD$ contains North-West quadrants 
$[\n', \8) \times (-\8, \b']$ with
$\n'>\n_c, \b'<\b_c^-((\n_c+\n')/2)$. 
Once (\ref{n_c<8}) is confirmed, (\ref{eq:betanuc}) follows 
immediately from the definition of $\n_c$ and the fact that 
$\cD$ is a closed set. 


\subsection{Critical curves}

We introduce:
\bdnl{alpha(b)}
\a (\b)={(e^\b-1)^2 \over e^\b (e^\b-1-\b)}, \; \; \; \b \in \R, 
\; \; \mbox{with $\a (0)\st{\rm def.}{=}2$.}
\edn
We note that
\bdnl{alpha'(b)}
\mbox{$\a (\b)$ decreases from $\a (-\8)=+\8$ to $\a (\8)=1$.}
\edn
Our first main result consists in upper and lower bounds on the critical values $\b^\pm_c(\n)$.
The following estimates show in particular that 
$\b^+_c (\n)$ (resp. $\b^-_c (\n)$) 
is locally Lipschitz continuous and strictly decreasing 
(resp. increasing) in $\n$. 
\Theorem{bc}
Let $d \ge 3$.
\bds
\item[(a1)] 
If $0<\n_0<\n$ and $1 \le \a \le \a (\b^+_c (\n_0))$, then, 
\bdnl{bc+sand+1}
\ln \lef( 1+c_1^+(\n_0)\lef( \n_0 \over \n\ri)^{1/\a}\ri) 
\le \b^+_c (\n) \le \ln \lef( 1+c_1^+(\n_0)\lef( \n_0 \over \n\ri)^{1/2}\ri),
\edn
\bdnl{bc+sand+2}
c_2^+(\n_0) \lef( 1-\lef( \n_0 \over \n\ri)^{1/2} \ri) 
\le \b^+_c (\n_0)-\b^+_c (\n) \le 
c_2^+(\n_0) \lef( \lef( \n \over \n_0\ri)^{1/\a}-1 \ri),
\edn
where 
$c_1^+(\n_0)=\lm ( \b^+_c (\n_0))$ and 
$c_2^+(\n_0)=1-\exp (-\b^+_c (\n_0))$. 
\item[(a2)] 
If $0<\n_1 \le \n<\n_0$ and $1 \le \a \le \a (\b^+_c (\n_1))$,  then,  
\bdnl{bc+sand-1}
\ln \lef( 1+c_1^+(\n_0)\lef( \n_0 \over \n\ri)^{1/2}\ri) 
\le \b^+_c (\n) \le \ln \lef( 1+c_1^+(\n_0)\lef( \n_0 \over \n\ri)^{1/\a}\ri),
\edn
\bdnl{bc+sand-2}
c_2^+(\n_0)\lef( 1-\lef( \n \over \n_0 \ri)^{1/2}  \ri) 
\le \b^+_c (\n)-\b^+_c (\n_0) \le 
c_2^+(\n_0) \lef( \lef( \n_0 \over \n\ri)^{1/\a}-1 \ri).
\edn
\item[(b1)] 
If $\n_c<\n_0<\n$ (cf. (\ref{nuc})) and $\a (\b^-_c (\n_0)) \le \a $, then,  
\bdnl{bc-sand+1}
\ln \lef( 1-c_1^-(\n_0)\lef( \n_0 \over \n\ri)^{1/\a}\ri) 
\le \b^-_c (\n) \le \ln \lef( 1-c_1^-(\n_0)\lef( \n_0 \over \n\ri)^{1/2}\ri),
\edn
\bdnl{bc-sand+2}
c_2^-(\n_0) \lef( 1-\lef( \n_0 \over \n\ri)^{1/\a} \ri)
\le \b^-_c (\n)-\b^-_c (\n_0) \le 
c_2^-(\n_0)\lef( \lef( \n \over \n_0 \ri)^{1/2}-1 \ri),
\edn
where 
$c_1^-(\n_0)=|\lm ( \b^-_c (\n_0))|$ and 
$c_2^-(\n_0)=\exp (-\b^-_c (\n_0))-1$. 
\item[(b2)] 
If $\n_c <\n_1\st{\rm def}{=}\n_0c_2^-(\n_0)^2<\n <\n_0$ 
and $\a \ge \a (\b^-_c (\n_1))$, then,  
\bdnl{bc-sand-1}
\ln \lef( 1-c_1^-(\n_0)\lef( \n_0 \over \n\ri)^{1/2}\ri) 
\le \b^-_c (\n) \le \ln \lef( 1-c_1^-(\n_0)\lef( \n_0 \over \n\ri)^{1/\a}\ri),
\edn
\bdnl{bc-sand-2}
c_2^-(\n_0) \lef( \lef( \n_0 \over \n\ri)^{1/\a}-1 \ri)
\le \b^-_c (\n_0)-\b^-_c (\n) \le 
c_2^-(\n_0)\lef( 1-\lef( \n \over \n_0\ri)^{1/2} \ri).
\edn
\eds
\end{theorem}
The proof of \Thm{bc} will be presented in section \ref{auxiliary}. 
From the above estimates we derive the following 
quantitative
informations. The first one  follows from (\ref{bc+sand-1}) 
with $\n=\n_1$ and $\a=1$, and the second one is from \Cor{decay*of*b_c}:
\Corollary{b_casymp} 
For $d \geq 3$, we have 
 \bdnl{bc(0)+}
\b^+_c (d, \n) \asymp \ln (1/\n)
\; \; \mbox{as $\n \searrow 0$},
\edn
and also,
\begin{equation} \label{eq:decay*of*b_c2}
|\b^\pm_c (d,\n)| \asymp 1/\sqrt{\n} \quad {\rm as } \quad \n \nearrow \8.
\end{equation}
\end{corollary}

\begin{remark} \label{rk:order}
In dimension $d \geq 3$, a phase transition occurs on the curve $\Crit$.
The function $p(\b,\n)$, being equal to the analytic function $\n \lambda(\b)$ on one side  of the curve,
takes a different value on the other side.  We do not know what is the order of the phase transition.
However, we will show that, when  $d \geq 3$, the gradient of $p_t(\beta,\nu)$ converges as $t \to \8$ to that of the
 limit $\n \lm(\b)$, at all  points of $\cD$, in particular those in  $\Crit$ 
(See the argument at the end of section \ref{p*th:phasediag}).
\end{remark}
\subsection{Path localization}
For all $t >0$ we define, in  \Prop{prop:measurable}
below, a  measurable function 
$(s,\h) \mapsto \Yt(s)$ with values in $\R^d$ such that
\begin{equation} \label{eq:Y2}
\mu_t\big(B_s \in U(
\Yt(s))\big)=
\max_{x \in \R^d} \mu_t\big(B_s \in U(x)\big), \qquad s\in [0,t].
\end{equation}
$\Yt$ depends on the environment $\h$, it is not continuous in general, 
however we call it the "optimal path" or "favorite path".
It is convenient to introduce the notation
$$
\chi_{s,x}={\mathbf 1}_{B_s\in U(x)}
$$
for the indicator function that the path sees the point $(s,x)$,
so that $\chi_{s,x}=\chi_{s,x}(B)={\mathbf 1}_{V_t}(s,x)$, and
$$
\h (V_t)= \int \chi_{s,x} \h_t(ds,dx) .
$$
We define  the overlap between two replicas $B$ and $\tilde B$, 
which plays a major role in quantitative estimates for localization:
\bdnl{R_t}
R_t=R_t(B, \tilde B)=\frac{1}{t}| V_t(B) \cap V_t(\tilde B)|
\edn
and the overlap between a polymer path and the optimal path,
\bdnl{R*(t)}
R^*_t=R^*_t(B, \h)=\frac{1}{t} \int_0^t  {\mathbf 1}_{\{B_s \in U(\Yt(s))\}} ds.
\edn
Note that $R, R^*$ take values in the interval $[0,1]$,  that
\bdnl{ER_t}
 \mu_t^{\otimes 2} (R_t)
=\frac{1}{t} \int_0^t ds\int_{\R^d} \mu_t( \chi_{s,x})^2 dx,
\edn
by Fubini's theorem, and
\bdnl{ER*(t)}
\mu_t (R^*_t)
=\frac{1}{t} \int_0^t  \mu_t\lef[  \chi_{s,
{\mathsf Y}^{(t)}(s)}\ri] ds
= \frac{1}{t} \int_0^t  \max_x \mu_t \lef[  
B_s \in U(x)
\ri] ds.
\edn

We say that $\b>0$ [resp., $\b<0$]
 is a {\it point of increase} of $\n \lm -p$ if
$\n \lm(\b') -p(\b',\n)>\n \lm(\b) -p(\b,\n)$ 
for all $\b'>\b$ [resp., $\b' < \b$].
If $\beta >0$ is a  point of increase 
for $\n \lm -p$, then necessarily $\b \geq \b^+_c$.
By monotonicity in \Thm{th:psi} (b), we already know that 
\bdnl{eq:conj}
\lef({\6 p \over \6 \b}\ri)_+ (\b,\n) 
\leq \n \lm'(\b)\;\;\; {\rm for \ }  \b \geq 0, \qquad
\lef({\6 p \over \6 \b}\ri)_{-} (\b,\n) 
\geq \n \lm'(\b) \;\;\;  {\rm for \ }  \b \leq 0,
\edn
where 
$\left(\frac{\6}{\6 \b}\right)_{\!\!-}$ [resp., 
$\left(\frac{\6}{\6 \b}\right)_{\!\!+}$] denote the left 
[resp., right] derivative. For a fixed $\nu$, 
$\beta \mapsto p(\beta, \nu)$ is differentiable except for at most 
countably many $\b$'s, and hence, 
$\left(\frac{\6 p}{\6 \b}\right)_{\pm}(\b,\n)
=\frac{\6 p}{\6 \b}(\b,\n)$ except such 
$\b$'s. 

The following result, in particular (\ref{eq:pathloc2}) which is the punchline, shows that a localization properties 
of the polymer is equivalent to the strictness of the 
inequality (\ref{eq:conj}), up to the exceptional 
non-differentiability of $p (\cdot, \nu)$.
\Theorem{th:pathloc}
\bds
\item[(a)] There exists $c=c(d) \in (0,1]$ such that
\bdnl{eq:pathloc3}
c\lef( Q \mu_t \left( R^*_t \right) \ri)^2 
\le Q \mu_t^{\otimes 2} \left( R_t \right) 
\le Q \mu_t \left( R^*_t \right).
\edn
\item[(b)]
For $\b \neq 0$ and $\n >0$, define 
\bdnl{del(bn)}
\del_\pm (\b,\n)=
(\n\lm)^{-1}  \left( \n \lm'(\b)-  
\lef(\frac{\6 p}{\6 \b} \ri)_{\pm}(\beta,\nu) \right) \ge 0\;,
\; \; \; 
\edn
(cf. (\ref{eq:conj})). Then, 
\bdnl{eq:pathloc4}
\liminf_{t\to\8} Q \mu_t^{\otimes 2} \left( R_t \right)
 \ge 
\lef\{ \barray{ll}
e^{-\b}\del_+ (\b,\n)
& \mbox{if $\b >0$}, \\
\del_- (\b,\n)
& \mbox{if $\b <0$}. 
\earray \ri. 
\edn
In particular, 
\bdnl{eq:pathloc2}
\liminf_{t\to\8} Q \mu_t \left( R^*_t \right) 
\ge  \liminf_{t\to\8} Q \mu_t^{\otimes 2} \left( R_t \right) 
>0
\edn 
whenever the inequality (\ref{eq:conj}) is strict. Moreover, 
\bdnl{eq:pathloc5}
\limsup_{t\to\8} Q \mu_t^{\otimes 2} \left( R_t \right)
 \le 
\lef\{ \barray{ll}
\del_- (\b,\n)
& \mbox{if $\b >0$}, \\
e^{|\b|}\del_+ (\b,\n)
& \mbox{if $\b <0$}. 
\earray \ri. 
\edn
\item[(c1)]
For a fixed $\nu$, $\del_-(\b, \nu) \ge \del_+(\b, \nu)>0$ if $\b>0$ is 
large enough. 
\item[(c2)] For all  points $\b$  of increase of $\n \lm -p (\cdot, \nu)$, 
there exists a sequence $(\b_k)_k$ converging to $\b$
such that for all $k$, 
$p (\cdot, \nu)$ is differentiable at $\b_k$ and 
$\del_+ (\b_k, \nu)=\del_- (\b_k, \nu)>0$.
\eds
\end{theorem}
The proof of \Thm{th:pathloc} will be presented in section 
\ref{p*pathloc}. 
Such statements express the strong localization properties of the polymer. 
In particular, the time-average $(1/t) \int_0^t  
{\mathbf 1}_{B_s\in U(\Yt(s))}ds$ is the {\em time fraction} the polymer spends with the favourite path. Under the strictness of the inequality (\ref{eq:conj}), 
the time fraction is positive. 
For a benchmark, we recall that, for the free measure $P$, 
for all smooth path ${\mathsf Y}$ and all $\delta>0$, there exists 
a positive $C$ such that for large $t$,
\bdnl{bench}
P\left( \frac{1}{t} \int_0^t  
{\mathbf 1}_{B_s\in U( {\mathsf Y}(s))}ds  \geq \delta \right)  \leq \exp \{-Ct\}
\edn
(In fact, it is not difficult to see (\ref{bench}) for $Y \equiv 0$ by 
applying Donsker-Varadhan's large deviations \cite{DoVa75}. 
Then, one can 
use Girsanov transformation to extend (\ref{bench}) to the case of 
smooth path ${\mathsf Y}$.)

The results of \Thm{th:pathloc} are to be compared 
with the results in \cite{CY05}, that we recall now. These ones only deal with end points,
i.e., with the location of the polymer at the last moment it interacts with the medium. 
If $p(\n,\b)<\n \lm(\b)$, we have by Theorem 
2.3.2 and Remark 2.3.1 of \cite{CY05} that 
\bdnl{eq:pathloc1}
\liminf_{t \to\8} \frac{1}{t} \int_0^t  
\max_{x \in \R^d} \mu_s( B(s) \in U(x))
ds 
 \ge 
\liminf_{t\to\8} \frac{1}{t} \int_0^t\int_{\R^d} \mu_s( \chi_{s,x})^2 ds dx >0,
\; \; \; \mbox{$Q$-a.s.}
\edn
Note that $\Yt(t)$ is a maximizer of $ \mu_t( B(t) \in U(x))$. 
As in \cite{CC11}, we conjecture that the set where localization occurs, coincides with the full
low temperature/high density region $\cL$, where the quenched free energy is strictly smaller than the annealed one:
\begin{conjecture}
$$\Big\{ (\n,\b): \mbox{The inequality (\ref{eq:conj}) is strict} \Big\}=\cL.$$
\end{conjecture}
We now turn to {\it complete path localization}. In the region where $\n \b^2$ is large with $\b$ bounded, 
the localization becomes strong in various aspects. First of all, 
the fraction of time 
$t^{-1} \int_0^t   {\mathbf 1}_{B_s\in U(\Yt(s))}ds $ that the polymer spends in the neighborhood 
of the favourite path 
  tends to its maximum value 1. 
This behavior is in  
a sharp contrast with the benchmark mentioned above.
\Theorem{th:pathlocfull}(Complete localization 1)
Let $\b_0 \in (0,\8)$ be arbitrary. 
Then, as $|\b| \le \b_0$ and $\n \b^2 \ra \8$, 
\bdmn
\liminf_{t\to\8} Q \mu_t \left( R^*_t \right) 
& \ge &
\liminf_{t\to\8} Q \mu_t^{\otimes 2} \left( R_t \right) 
\label{eq:pathlocfull2}\\
 & = & 1-\cO \lef( (\n \b^2)^{-1/6}\ri).
\label{eq:pathlocfull1}
\edmn
\end{theorem}
In other words, for any bounded function $\b(\nu)$ with $\n \b(\n)^2 \to \8$ as $\n \to \8$, 
the limit 
$\ell(\b,\n)= \liminf_{t\to\8} Q \mu_t \left( R^*_t \right)  \in [0,1]$
converges to its maximal value,
$$
\ell (\b(\n),\n) \to 1, \quad \n \to \8.
$$
Of course, the parameter $\b$ can become small, but not too much, since no localization occurs if $\b =0$.
The proof of \Thm{th:pathlocfull}, together with that of 
\Thm{th:pathlocfull2} below, will be presented in section 
\ref{p*comploc}. 
\medskip

We next 
extract fine additional  information on the 
geometric properties of the Gibbs measure.
For $\delta \in (0,1/2)$ define the $(\delta,t)$-negligible set as
$$
{\cal N}_{\delta,t}^\h =  \Big\{(s,x) \in [0,t]\times \rd  : \mu_t(\chi_{s,x}) \leq \delta\Big\}, 
$$ 
and the  $(\delta,t)$-predominant set as
$$
{\cal P}_{\delta,t}^\h =  \Big\{(s,x) \in [0,t]\times \rd  : \mu_t(\chi_{s,x}) \geq 1-\delta\Big\}.
$$ 
As suggested by the names, ${\cal N}_{\delta,t}^\h$ is the set of space-time locations the polymer wants to stay
away from, and ${\cal P}_{\delta,t}^\h$ is  the set of  locations the polymer likes to visit. Both sets depend on the environment.
\Theorem{th:pathlocfull2}(Complete localization 2)
For all $0<\delta <1/2$, we have as  $|\b| \le \b_0$ and $\n \b^2 \ra \8$, 
\begin{equation}
\label{eq:r2to1-3}
\limsup_{t\to\8}Q \frac{1}{t}  \Big\vert  ({\cal N}_{\delta,t}^\h \cup  {\cal P}_{\delta,t}^\h)^\complement
 \Big\vert  
=\cO \lef( (\n \b^2)^{-1/6}\ri),
\end{equation}
\begin{equation}
\label{eq:r2to1-4}
\limsup_{t\to\8}Q \mu_t \left[ \frac{1}{t}  \Big\vert  V_t(B) \bigcap {\cal N}_{\delta,t}^\h
\Big\vert  
\right]
=\cO \lef( (\n \b^2)^{-1/6}\ri),
\end{equation}
\begin{equation}
\label{eq:r2to1-5}
\limsup_{t\to\8}Q \mu_t \left[ \frac{1}{t}  \Big\vert  
V_t(B)^\complement   \bigcap 
{\cal P}_{\delta,t}^\h
  \Big\vert  
\right]
=\cO \lef( (\n \b^2)^{-1/6}\ri).
\end{equation}
\end{theorem}
Recall that $\vert \cdot \vert$ denotes the Lebesgue measure on $\R_+ \times \R^d$, and  note that 
$ \vert  {\cal N}_{\delta,t}^\h  \vert =  \vert  V_t(B)^\complement  \vert = \8$.

The limits (\ref{eq:r2to1-3}),  (\ref{eq:r2to1-4}),  (\ref{eq:r2to1-5}), bring 
information on how is the corridor around the favourite path where the measure concentrates for large $\nu \b^2$.
In this limit:
\begin{itemize}
\item
 most (in Lebesgue measure) time-space  locations become negligible or predominant, 
\item most  (in Lebesgue and Gibbs measures)  negligible  locations are outside the tube around the polymer path, 
\item
 most  (in Lebesgue and Gibbs measures)  predominant  locations are inside the tube around the polymer path.
\end{itemize}
The trace $ \big\{x \in \R^d: \mu_t(\chi_{s,x}) \geq 1-\delta\big\}$ at time $t$ of the $(\delta,t)$-predominant set is reminiscent of the $\epsilon$-atoms discovered in the discrete setting in \cite{Va07}, with $\epsilon=1-\delta$.
%
%
\SSC{Replica overlaps and favourite path}
We build on ideas similar to \cite{CC11}. Since the state space is continuous, some measurability issues appear, but also the geometric properties of the path measure are of interest.
\subsection{Favourite path}
For all times $s \leq t$, the function  
$x \mapsto \mu_t(B_s \in U(x))=\mu_t(\chi_{s,x})$ achieves 
its maximum, and the set of  maximizers is compact. 
We want to consider ''the maximizer", 
by  selecting a specific element in the 
argmax in case of multiplicity, but this
can be effective only with some measurability property. 
For a function $f: \R^d \to \R$ and a set $A \sub \rd$, we denote by
\bdnl{argmax}
\arg \max_{x \in A} f(x) = \{ x \in A: f(x)=\sup_{z \in A} f(z)\}
\edn
the set of the maximizer of $f$ on $A$.
\Proposition{prop:measurable}
There exists a measurable subset $\cM^0 \subset \cM$, 
and for each fixed $t >0$, 
a measurable function
$$
(s,\h) \mapsto \Yt(s): [0,t] \times \cM^0 \to \R^d
$$ 
such that
\bdmn
& & \mbox{$Q(\cM^0)=1$ and $Z_t (\h)<\8$ for all 
$\b \in \R$, $t>0$ and $\h \in \cM^0$,} 
\label{Z_t<8} \\
& & \Yt(s)\in  \arg \max_{x \in \R^d} \mu_t\big(B_s \in U(x)\big).
\label{eq:Y}
\edmn
\end{proposition}
As indicated in the notation, 
we will regard $\Yt(\cdot)$ as a function on $[0,t]$, 
which depends on $t$, but we keep in mind that it also depends 
on 
$\b$ and on $\h \in \cM^0$. It is not continuous in general, 
however we call it the "favourite path". 

\vvs
\noindent Proof of \Prop{prop:measurable}: 
We recall that 
$\cM$ is a Polish space  
under the vague topology ${\cal T}$ 
\cite[p. 170, 15.7.7]{Kal83} and that the Borel $\s$-field 
$\s [{\cal T}]$ coincides with $\cG$ \cite[p. 32, Lemma 4.1]{Kal83}.
This observation enables us to exploit a measurable selection theorem 
from \cite[p.289, Theorem 12.1.10]{StVa79}, as we explain now. 

We start with the definition of $\cM^0$. Let $r=r_d$ be the radius of $U(0)$ and 
$$
\tl{V}_t=\{ (s,x) \in [0,t] \times \rd \; ; \; 
  \exists u \geq 0, |u-s| \leq 1 ,  |x-B_u| \le r+1\}
$$
be an enlargement of $V_t$. 
We define 
\bdnl{cM0}
\cM^0 =\bigcap_{\b ,t >0}\left\{\h \in \cM\; ; \; 
P\big[\exp (\b \h (\tl{V}_t))\big]<\8\right\},
\edn
which satisfies (\ref{Z_t<8}). 
In the 
following argument, we always assume that $\h \in \cM^0$. 
 With 
$$
Z_t(s,x,\h)=P[ \z_t(\h); B_s \in U(x)],\; \; \;  \z_t(\h)=\exp\{\b \h(V_t)\},
$$ 
we can write the  right-hand side of (\ref{eq:Y}) as
$$
 \arg \max_{x \in \R^d} \mu_t\big(B_s \in U(x)\big)
= \arg \max_{x \in \R^d} Z_t(s,x,\h).
$$
Thus, 
we wish to select a maximizer 
of $x \mapsto Z_t(s,x,\h)$ as a 
measurable function in $(s,\h)$. 
As can be seen below, our method relies heavily on the continuity 
of the functions we work with (cf. the proof of \Lem{K(sh)}). 
Unfortunately, 
$Z_t(s,x,\h)$ is discontinuous at all $(s,x,\h)$ 
such that $\h(\{s\}\times U(x))\geq 1$. 
To circumvent 
this obstacle, we will consider a continuous approximation 
of $Z_t(s,x,\h)$ (cf. (\ref{Z^(k)}) below), 
together with a cut-off of $x$-variables. 
Let $\varphi: \R \to \R_+$ be continuous, 
supported inside $[-1,1]$, and the integral equal to one. 
For $k\geq 1$, let 
$\vak(t)=k\varphi(kt)$, and $\pk(x)=\big[1-k\times {\rm dist}(x, U(0))\big]^+$.
Then we have, as $k \to \8$, $\vak(t) dt \to \del_0(dt)$ 
weakly, and $\pk \to  {\mathbf 1}_{U(0)}$ pointwise.
Define 
$$
\zk_t(\h)=\exp \{ \b \int \h_t(dsdx) \int_\R \vak (s-u) \pk (B_{u}-x)du\},
$$
and
\bdnl{Z^(k)}
Z^{(k)}_t(s,x,\h)=P[ \zk_t(\h); B_s \in U(x)]. 
\edn
Note that, with $u^+=\max\{u,0\}$,  
\bdnl{zkt<}
\zk_t(\h) \le \exp (\b^+\h (\tl{V}_{t+1})),
\edn 
and hence 
$Z^{(k)}_t(s,x,\h) <\8$ for all $\h \in \cM^0$.  
Let $\cK$ be the totality of compact subsets in $\rd$, equipped with the 
Hausdorff metric. 
Then, we will show in \Lem{K(sh)} below that,  for any integer $\ell \geq 1$, the 
mapping,    
\bdnl{K(sh)}
(s,\h) \mapsto K^{(k, \ell)}(s,\h) = 
\arg \max_{x \in [-\ell, \ell]^d}  Z^{(k)}_t(s,x,\h)
\edn 
defined by (\ref{argmax}), 
is Borel measurable from $[0,t] \times \cM^0$ to $\cK$. 
Thanks to this measurability, which we will assume for the moment, 
we deduce from the measurable selection theorem mentioned 
above, that there exists a measurable mapping 
${\mathsf Y}^{(t, k, \ell)}: [0,t]\times \cM^0
\to  [-\ell,\ell]^d$ such that 
$$
{\mathsf Y}^{(t, k, \ell)}(s,\h) \in K^{(k, \ell)}(s,\h). 
$$
We now let $k \ra \8$. First, we see from a standard mollifier argument 
that $\zk_t (\h)\st{k \to \8}{\to} \z_t (\h)$ 
for fixed $\h$ and $B$. Thus, we have by
(\ref{zkt<}) and the dominated convergence theorem that 
\begin{equation} \label{eq:2v}
\forall \h \in \cM^0, 
\qquad Z^{(k)}_t(s,x, \h) \st{k \ra \8}{\to} Z_t(s,x,\h) 
\quad \mbox{uniformly in $(s,x) \in [0,t] \times \rd$}.
\end{equation}
This implies that every
limit point of ${\mathsf Y}^{(t, k, \ell)}(s,\h) $ as $k \to \8$ 
is a maximizer of $x \ra Z_t(s,x,\h)$ on $[-\ell,\ell]^d$. 
We construct such a limit point in a measurable manner: 
the lower limit of the first component $\mathfrak L_1=\liminf_k {\mathsf Y}_1^{(t, k, \ell)}$
is measurable in $(s,\h)$, and we can define an  extractor  by
$$k(m)=\inf\{k >  k(m-1):  {\mathsf Y}_1^{(t, k, \ell)}\leq \mathfrak L_1 + (1/m)\},$$ which depends on $(s,\h)$ in a 
measurable manner. We now restrict to the 
sequence $({\mathsf Y}^{(t, k(m), \ell)}; m\geq 1)$ of measurable functions, whose first coordinate converges.
Now we repeat extraction for the other coordinates, starting with the second one. We end up with a converging subsequence that we still
denote by the same symbol   $({\mathsf Y}^{(t, k(m), \ell)}; m\geq 1)$ which  converges  pointwise  as $m \to \8$ to some
 ${\mathsf Y}^{(t, \ell)}$, a maximizer of  $Z_t(s,x)$ on $[-\ell,\ell]^d$, and which is measurable. 
Note now that 
$$
\ell(s,\h)= \inf\{\ell \geq 1: \max_{x \in [-\ell,\ell]^d} Z_t(s,x)
= \max_{x \in \R^d} Z_t(s,x)\},
$$
is measurable. Hence, it suffices to take 
$\Yt(s)={\mathsf Y}^{(t, \ell(s,\h))}$, this ends the proof.
\hfill $\Box$ 
\medskip
\Lemma{K(sh)}
The mapping (\ref{K(sh)}) is Borel measurable. 
\end{lemma}
Proof: 
We approximate the set $\cM^0$ (cf. (\ref{cM0})) by:
$$
\cM^0_{\b,t,L} \st{\rm def.}{=}\{\h \in \cM^0\; ; \; 
P\exp (2\b^+ \h (\tl{V}_{t+1}))  \le L\} \nearrow \cM^0, 
\; \; \; L \nearrow \8.
$$
Here, the parameters $\b$ and $t$ of $\cM^0_{\b,t,L}$ 
are the same as those of $Z^{(k)}_t(s,x,\h)$. 
It is now enough to prove the Borel measurability of $K^{(k, \ell)}(s,\h)$ 
on $[0,t] \times \cM^0_{\b,t,L}$ for any $L$. 
For such measurability, the following sufficient condition is 
known, cf. \cite[p.289, Lemma 12.1.8]{StVa79}:
\bdnl{lemSV}
\barray{l}
\mbox{For any sequence 
$(s_n,\h_n) \ra (s_0,\h_0)$ in $[0,t] \times \cM^0_{\b,t,L}$ 
and $x_n \in K^{(k, \ell)}(s_n,\h_n)$, it is} \\
\mbox{true that $x_n$ has a limit point $x_0$ in $K^{(k, \ell)}(s_0,\h_0)$.}
\earray
\edn
Let us verify the above criterion. We start by noting that, for fixed $k$,
\begin{equation} \label{eq:1v}
 \mbox{$(s,x,\h) \ra Z^{(k)}_t(s,x,\h)$ is continuous on 
$[0,t] \times \rd \times \cM^0_{\b,t,L}$}. 
\end{equation}
In fact, suppose that 
$(s_n,x_n,\h_n) \ra (s,x,\h)$ in $[0,t] \times \rd \times \cM^0_{\b,t,L}$. 
We write:
\bdnn
\lefteqn{|Z^{(k)}_t(s_n,x_n,\h_n)-Z^{(k)}_t(s,x,\h)|} \\ 
& \le & 
|Z^{(k)}_t(s_n,x_n,\h_n)-Z^{(k)}_t(s_n,x_n,\h)|
+ |Z^{(k)}_t(s_n,x_n,\h)-Z^{(k)}_t(s,x,\h)| \\
& \le & I_n+J_n,
\ednn
where
$$
I_n=P[|\zk_t(\h_n)-\zk_t(\h)|],
\; \; \; 
J_n =P[\zk_t(\h)^2]^{1/2}P[|\chi_{s,x}-\chi_{s_n,x_n}|^2]^{1/2}.
$$
We have that $I_n\st{n \ra \8}{\ra}0$, 
since $\zk_t(\h_n) \st{n \ra \8}{\ra}\zk_t(\h)$ for fixed $B$, 
and $\{\zk_t(\h)\; ; \; \h \in \cM^0_{\b,t,L}\}$ is uniformly integrable. 
On the other hand, we have $J_n\st{n \ra \8}{\ra}0$, since 
$P[\zk_t(\h)^2] <\8$ 
for $\h \in \cM^0_{\b,t,L}$, and $P(B_s \in \6 U (x))=0$. 

To verify (\ref{lemSV}), 
let $(s_n,\h_n)$, $(s_0,\h_0)$ and $x_n$ be as its assumption. 
Since $[-\ell, \ell]^d$ is compact, we can take a converging 
subsequence 
$x_{n(j)} {\ra} x_0$ as $j {\ra} \8$.
On the other hand, we see as a consequence of (\ref{eq:1v}) that
\bdnl{max*lsc}
\mbox{
${\dps (s,\h) \mapsto \max_{x \in [-\ell, \ell]^d}Z^{(k)}_t(s,x,\h)}$ 
is lower semi-continuous}.
\edn
Hence, 
\bdnn
Z^{(k)}_t(s_0, x_0, \h_0) & \st{\mbox{\scriptsize (\ref{eq:1v})}}{=} & 
\lim_{j \ra \8}Z^{(k)}_t(s_{n(j)}, x_{n(j)}, \h_{n(j)}) \\
& = & \lim_{j \ra \8}\max_{x \in [-\ell, \ell]^d}Z^{(k)}_t(s_{n(j)},x,\h_{n(j)}) \\
& \st{\mbox{\scriptsize (\ref{max*lsc})}}{\ge} & 
\max_{x \in [-\ell, \ell]^d}Z^{(k)}_t(s_0,x,\h_0).
\ednn
Thus, we have verified  (\ref{lemSV}). \hfill $\Box$ 
\subsection{Overlaps} 
In the next section, we will obtain information in a two-replica system, hence on the value of 
$R_t$. In order to translate it into one for a single path of the polymer measure, we will use an elementary lemma, where the first item takes care of positive overlaps, and the second one of values close to 1. Recall definitions (\ref{R_t}) and  (\ref{R*(t)}) of $R_t, R^*_t$.
\Lemma{lem:2to1}
{\bf (Two-to-one lemma)}
 Almost surely, we have the following: \\
{\rm (i)} $\exists c=c(d) \in (0,1)$ such that 
\begin{equation}
\label{eq:2to1-1}
c \;\mu_t( R^*_t)^2 \leq \mu_t^{\otimes 2} ( R_t) \leq \mu_t( R^*_t);
\end{equation}
{\rm (ii)}
\begin{equation}
\label{eq:2to1-2}
\mu_t\Big( 1-R^*_t\Big) \leq \mu_t^{\otimes 2} \Big(1- R_t\Big).
\end{equation}
Moreover, for all $\delta \in (0,1/2]$,
\begin{eqnarray}
\label{eq:2to1-3}
\frac{1}{t}  \Big\vert  \Big\{(s,x) \in [0,t]\times \rd  : \mu_t(\chi_{s,x}) \in [\delta, 1-\delta]\Big\} \Big\vert  &\leq&
\frac{1}{\delta (1-\delta)} \;
 \mu_t^{\otimes 2} \Big(1- R_t\Big),\\
 \label{eq:2to1-4}
\mu_t \left[ \frac{1}{t}  \Big\vert  V_t(B) \bigcap \Big\{(s,x) : \mu_t(\chi_{s,x}) \leq \delta\Big\} 
\Big\vert  
\right]  &\leq&
\frac{1}{1-\delta} \;
 \mu_t^{\otimes 2} \Big(1- R_t\Big),\\
\label{eq:2to1-5}
\mu_t \left[ \frac{1}{t}  \Big\vert  
V_t(B)^\complement   \bigcap 
 \Big\{(s,x) : \mu_t(\chi_{s,x}) \geq 1-\delta\Big\} 
  \Big\vert  
\right] &\leq&
\frac{1}{1-\delta} \;
 \mu_t^{\otimes 2} \Big(1- R_t\Big).
\end{eqnarray}
\end{lemma}

\vvs
\noindent Proof of the Lemma: 
Since
\begin{eqnarray*}
\mu_t^{\otimes 2} ( R_t) &=& \frac{1}{t} \int_0^t \int \mu_t( \chi_{s,x})^2  dx ds \\
&\leq& \frac{1}{t} \int_0^t ds  \max_x \mu_t \lef[  
B_s \in U(x)
\ri] \times
\int_{\rd} \mu_t( \chi_{s,x})  dx 
\\&=& \mu_t( R^*_t)
\;,
\end{eqnarray*}
that is the right-hand-side inequality of (\ref{eq:2to1-1}).
To prove the left-hand-side, we introduce a smaller ball
$\frac{1}{2}U(0)=\{ \half z \; ;\; z \in U(0)\}$.
By the Schwarz inequality,
\bdnn
\int_{\rd} \m_t \lef[B_s \in U(z)\rig]^2 dz
&\geq &
\lef| \half U(0) \rig|^{-1}
\lef(\int_{y+\half U(0)} \m_t \lef[B_s \in U(z)\rig] dz \rig)^2\\
&\geq & 2^d
\lef(\int_{y+\half U(0)} \m_t \lef[B_s \in y+\half U(0)
\rig] dz \rig)^2\\
&=& 2^{-d} \m_t \lef[B_s \in y+\half U(0)
\rig]^2\;,
\ednn
where we have used the triangular inequality in the second line. 
By additivity of $\m_t$,
$$
 \max_{y \in \rd} \m_t \lef[B_s \in U(y)\rig] \leq
c' \max_{y \in \rd} \m_t \lef[B_s \in y+\half U(0)\rig]\;,
$$
with $c'=c'(d)$ the minimal number of translates of  $\frac{1}{2}U(0)$
necessary to cover  $U(0)$. 
Combining these two estimates and integrating on $[0,t]$, we 
complete the proof of of (\ref{eq:2to1-1}).

The claim (\ref{eq:2to1-2}) is a reformulation of 
the second one in  (\ref{eq:2to1-1}).
The last claims follow from the inequality
$$
u(1-u) \geq (1-\delta)u {\mathbf 1}_{u < \delta} + 
\delta (1-\delta) {\mathbf 1}_{u \in[\delta,1-\delta]} +(1-\delta)(1-u) {\mathbf 1}_{u > 1-\delta} .
$$
Setting $A_s=\{x : \mu_t(\chi_{s,x} )\in [\delta, 1-\delta]\}$ and writing
\begin{eqnarray*}
\mu_t^{\otimes 2} (1- R_t) &=& \frac{1}{t} \int_0^t \int_{\rd} \left[ \mu_t( \chi_{s,x})- \mu_t( \chi_{s,x})^2\right] dx ds \\
&\geq &  \frac{1}{t} \int_0^t \int_{A_s} \left[ \mu_t( \chi_{s,x})- \mu_t( \chi_{s,x})^2\right] dx ds \\
&\geq & 
\delta (1-\delta) 
\frac{1}{t} \int_0^t \left| \{x : \mu_t(\chi_{s,x}) \in [\delta, 1-\delta]\} \right|
 ds, 
\end{eqnarray*}
which yields (\ref{eq:2to1-3}). For the next one, we write
\begin{eqnarray*}
\mu_t^{\otimes 2} (1- R_t) &=& \frac{1}{t} \int_0^t \int_{\rd} \left[ \mu_t( \chi_{s,x})- \mu_t( \chi_{s,x})^2\right] dx ds \\
&\geq& (1-\delta) \frac{1}{t} \int_0^t \int_{\rd}  \mu_t( \chi_{s,x})  {\mathbf 1}_{
 \mu_t( \chi_{s,x})  < \delta} ds dx\\
 &=& (1-\delta) \mu_t \left[ \frac{1}{t} \int_0^t \int_{\rd}   {\mathbf 1}_{
 \mu_t( \chi_{s,x})  < \delta, B_s \in U(x)}   \right] ds dx,
\end{eqnarray*}
which is (\ref{eq:2to1-4}). The last claim can be proved similarly. \hfill $\Box$ 
\medskip
\SSC{The arguments of the proof of path localization}
We need estimates on the free energy $p(\beta, \nu)$ and/or its derivative. By definition of the critical values, we have strict inequality 
between the quenched  free energy $p(\beta, \nu)$ and the annealed one $\nu \lambda(\b)$, which was enough
when dealing with fixed parameters such  that $\b \notin [\beta^-_c,\beta^+_c]$ to get end-point localization
via semi-martingale decomposition. For path localization, we use an integration by parts formula.
\subsection{Integration by parts}
By an elementary computation, we see that, for a Poisson variable $Y$ with parameter $\theta$, 
the identity ${\bf E} Yf(Y) = \theta {\bf E}
f(Y+1)$ holds for all non negative function $f$. This is in fact the integration by parts formula
for the Poisson distribution,  it implies  the first property below, already used in \cite{CY05},
that we complement by a second formula, more alike to usual integration by parts formulas.
\Proposition{ipp}
(i) For $h: [0,t]\times \rd \times \cM \to \R_+$ a 
measurable function, 
we have
\bdnl{ipp1}
Q \left[ \int h(s,x;\eta_t) \eta_t(dsdx) \right] =
\n \int_{[0,t]\times \rd} dsdx   
Q \left[ h(s,x; \eta_t + \delta_{s,x}) \right]\;.
\edn
(ii)  Let $h: [0,t]\times \rd \times \cM \to \R$ be a 
measurable function, such that there exists a compact $K \subset \R^d$ with 
$h(s,x;\h)=0$ for all $s \leq t, \h \in \cM, x \notin K$, and such that 
$\int_{[0,t]\times \rd} dsdx    Q \left[ |h(s,x; \eta_t)| \right] <\8$.
Then, with $\tl{\h}_t(dsdx)=\h (dsdx)-\n dsdx$, we have
\bdnl{ipp2}
Q \left[ \int  h(s,x;\eta_t) \tl{\h}_t(dsdx) \right] =
\n \int_{[0,t]\times \rd} dsdx   
Q \left[ h(s,x; \h_t + \delta_{s,x})-  h(s,x;\eta_t)\right]\;.
\edn
\end{proposition} 
Proof: Recall the (shifted) Palm measure $Q_{s,x}$ of 
 the point process $\h_t$, 
which can be thought of 
as the law of $\eta_t$ ``given that $\h_t\{(s,x)\}=1$'':
By definition of the Palm measure,
$$
Q [ \int h(s,x;\eta_t) \eta_t(dsdx) ] =
\n \int_{[0,t]\times \rd} dsdx \int_{\cal M} h(s,x; \h ) Q_{s,x}(d\h ) \;.
$$
By Slivnyak's theorem \cite[page 50]{SKM87} 
for the Poisson  point process $\h_t$, the Palm measure $Q_{s,x}$
is the law of $\h_t + \delta_{s,x}$, hence the  right-hand-side of the
above formula
is equal to the right-hand-side  of (i). 
The equality (ii) follows by considering the positive and
negative parts of $h$. \hfill $\Box$
Define
\bdnl{psi=lm-vp}
\qquad \widehat{p}_t(\b, \nu)
=p_t(\b, \nu)-\nu \b ,
\edn
which is a convex function of $\beta$.
By differentiation and using \Prop{ipp} we obtain the following 
\Lemma{energy}
For all $\beta \in \R$,
\bdnl{energy-eq}
t\;\frac{\6 p_t}{\6 \b} 
(\beta,\nu)  =\n  e^{\beta} \int_{[0,t]\times \rd} dsdx
\; Q \; \frac{\mu_t( \chi_{s,x}) }{1+ \lm \mu_t( \chi_{s,x})} \;,
\edn
and therefore 
\bdnl{3100}
t\;\frac{\6 \hat p_t}{\6 \b} 
(\beta,\nu)  =\n  \lm \int_{[0,t]\times \rd} dsdx
\; Q \; \frac{\mu_t( \chi_{s,x})-\mu_t( \chi_{s,x})^2 }{1+ \lm \mu_t( \chi_{s,x})} 
\;,\edn
\bdnl{3101}
t\;\frac{\6}{\6 \b} \big(\nu \lm(\b)- p_t (\beta,\nu) \big)
=\n e^\b  \lm \int_{[0,t]\times \rd} dsdx
\; Q \; \frac{[\mu_t( \chi_{s,x})]^2 }{1+ \lm \mu_t( \chi_{s,x})} 
\;.\edn
\end{lemma}
Proof: With the identity
$({\6}/{\6 \b}) \ln Z_t
=\mu_t( \h(V_t) )$, we obtain
from Fubini's theorem and \Prop{ipp},
\bdmn
t \;\frac{\6 p_t}{\6 \b} 
(\beta,\nu)  &=&
Q[ \mu_t( \h(V_t) )] \nn \\&=&
Q \int \h_t(dsdx) \mu_t[ \chi_{s,x}]\nn \\
&=& 
Q \int \h_t(dsdx) \frac{P[ \chi_{s,x}e^{\beta \h(V_t)}]}
{P[ e^{\beta \h(V_t)}]}\nn \\
&\st{\mbox{\scriptsize (\ref{ipp1})}}{=}& 
\n Q \int_{[0,t] \times \rd} dsdx  
\frac{P [\chi_{s,x}e^{\beta (\h+\delta_{s,x})(V_t)} ] }
{P [e^{\beta 
(\h+\delta_{s,x})(V_t)}]}\nn \\
&=& 
\n Q \int_{[0,t] \times \rd} dsdx   
\frac{e^{\beta} P [ \chi_{s,x}e^{\beta \h(V_t)}]}
{P[ (\lm\chi_{s,x}+1)e^{\beta \h(V_t)}]}\nn \\
&=& 
\n e^{\beta} 
Q \int_{[0,t] \times \rd} dsdx   \frac{ \mu_t
 [\chi_{s,x}]}
{1 + \lm \mu_t
 [\chi_{s,x}]}\;, \label{6p/6b}
\edmn
which is the first claim. 
Since $t=\int_{[0,t] \times \rd} dsdx  \mu_t
 [\chi_{s,x}]$,  we can express the left-hand side of  (\ref{3100}) 
as $\nu Q \int \int ds dx \psi(\mu_t(\chi_{s,x}))$, where
\bdnl{eq:psif}
\psi( u) = e^\b \frac{u}{1+\lambda u}-u=  \lambda \frac{u-u^2}{1+\lambda u}
\edn
by definition of $\lambda$, and similarly,  
the left-hand side of  (\ref{3101}) 
as $\nu Q \int \int ds dx \phi(\mu_t(\chi_{s,x}))$, where
\bdnl{eq:phif}
\phi( u) = -e^\b \frac{u}{1+\lambda u}
+e^\b u= e^\b \lambda \frac{u^2}{1+\lambda u}
\edn
\hfill $\Box$

For further use, note that we have for $u \in [0,1]$ and $\b \in \R$,
\bdmn
e^{-\b^+} \lambda (u-u^2)\leq \psi(u) & \leq & e^{-\b^-}\lambda (u-u^2), 
\label{eq:psif>} \\
\lm u^2 \le \phi(u) & \leq & e^\b \lm u^2. \label{eq:phif<}
\edmn
\subsection{Proof of path localization,  \Thm{th:pathloc} }
\label{p*pathloc}
From \Lem{energy} we can easily recover the following inequalities.
\begin{lemma} 
We have
$$
\n \leq \left(\frac{\6 p}{\6 \b}\right)_{\!\!-} (\b,\n)\leq   
\left(\frac{\6 p}{\6 \b}\right)_{\!\!+} (\b,\n)
\leq  \n \lm'(\b), \qquad 
\b \geq 0,
$$
and
$$
\n \geq 
\left(\frac{\6 p}{\6 \b}\right)_{\!\!+} (\b,\n)
\geq   
\left(\frac{\6 p}{\6 \b}\right)_{\!\!-} (\b,\n)
\geq  \n \lm'(\b), \qquad 
\b \leq 0.
$$
\end{lemma}
Proof: 
Inequalities follow from (\ref{3100}),  (\ref{3101}) 
and the convexity. $\Box$
\medskip

We now turn to the:\\
Proof of \Thm{th:pathloc}:
We consider the case $\b > 0$, the other case being similar.\\
(a) follows from (\ref{eq:2to1-1}) and Jensen's inequality. \\
(b): We first note that 
$$
e^{-\b}Q\m_t^{\otimes 2}(R_t)
\le
\frac{1}{t}\;\int_{[0,t]\times \rd} dsdx
\; Q \; \frac{\mu_t( \chi_{s,x})^2 }{1+ \lm \mu_t( \chi_{s,x})}
\le
Q\m_t^{\otimes 2}(R_t)
$$
by (\ref{ER_t}) and (\ref{eq:phif<}). 
 On the other hand, 
let $f_n:\R \ra \R$, $n=1,2,...$ be a sequence of 
convex functions which converges to a function $f$ pointwise. 
Then, it is easy to see that
\bdnl{lim*Dconv}
\lef( {d f\over d \b}\ri)_{-} \le 
\liminf_{n \ra \8}\lef( {d f_n\over d \b}\ri)_{-} \le 
\limsup_{n \ra \8}\lef( {d f_n\over d \b}\ri)_{+} \le 
\lef( {d f\over d \b}\ri)_{+}.
\edn
This, together with (\ref{3101}),  can be used as follows:
\bdnn
\n \lm e^\b\liminf_{t \ra \8} \frac{1}{t}\;\int_{[0,t]\times \rd} dsdx
\; Q \; \frac{\mu_t( \chi_{s,x})^2 }{1+ \lm \mu_t( \chi_{s,x})} 
&\st{\mbox{\scriptsize (\ref{3101})}}{=}&  \n \lm'(\b)
- \limsup_{t \ra \8} \frac{\6 p_t}{\6 \b} (\beta,\nu)  \\
&\st{\mbox{\scriptsize (\ref{lim*Dconv})}}{\ge} &
 \n \lm'(\b)-  \lef(\frac{\6 p}{\6 \b} \ri)_{+}(\beta,\nu), \\
\n \lm e^\b\limsup_{t \ra \8} \frac{1}{t}\;\int_{[0,t]\times \rd} dsdx
\; Q \; \frac{\mu_t( \chi_{s,x})^2 }{1+ \lm \mu_t( \chi_{s,x})} 
&\st{\mbox{\scriptsize (\ref{3101})}}{=}&  \n \lm'(\b)
- \liminf_{t \ra \8} \frac{\6 p_t}{\6 \b} (\beta,\nu)  \\
&\st{\mbox{\scriptsize (\ref{lim*Dconv})}}{\le} &
 \n \lm'(\b)-  \lef(\frac{\6 p}{\6 \b} \ri)_{-}(\beta,\nu)
\ednn
Putting things together, we get (\ref{eq:pathloc4}) and (\ref{eq:pathloc5}).\\
(c1):We see from the proof of  
\cite[(2.9)]{CY05} that, for a fixed $\n >0$, and $\b >0$ large,
\bdnl{cy05*1}
p (\b, \n) \le C_1\sqrt{\n \lm },\; \; \; C_1=C_1 (d) \in (0,\8).
\edn 
Suppose that $\b >0$ is sufficiently large. 
Then, by the convexity of $p (\cdot,\n)$ and (\ref{cy05*1}), 
\bdnn
\lef(\6 p \over \6 \b\ri)_{+}(\b,\n) 
& \le &  
 p (\b+1,\n)-p (\b,\n) 
 \\&\le&  C_2\sqrt{\n \lm },\; \; \; \qquad (C_2=C_2 (d) \in (0,\8)) \\
& \le & \n e^{\b},
\ednn
hence $\del_+(\b,\n)>0$. \\
(c2): By convexity, $p(\cdot,\n)$ 
is almost everywhere differentiable.
If $\b$ is a point of increase, we have for $\b_0 >\b$, 
$$
0< \n \lambda(\b_0) - p(\b_0,\n)-\n\lambda(\b)  + p(\b,\n) = 
\int_{\b}^{\b_0} \left( \n \lambda'(b)-\frac{\6 p}{\6 \b} (b,\nu) \right) db.
$$
Then the set of $b \in (\b,\b_0)$ such that the integrand is positive has non zero Lebesgue measure.
Since this holds for all $\b_0 > \b$, we conclude that there exists a decreasing sequence $\b_k \searrow \b$ with the desired properties. This ends the proof.
  \hfill $\Box$
\subsection{Proof of complete localization}
\label{p*comploc}
Complete localization holds in the asymptotics 
$|\b| \le \b_0 <\8, \n\b^2 \to \8$ because the quenched free 
energy diverges in a slower manner 
than the annealed one. Precisely, we will establish 
the following asymptotic estimate, which is key for a number of our results.
\Theorem{boundonp} 
Let $\b_0 \in (0,\8)$ be arbitrary. Then, 
\begin{equation}\label{eq:calO(...)}
 p(\nu, \beta ) = \beta \nu + \cO \big( (\nu \b^2)^{5/6}\big),
\; \; \mbox{as $|\b| \le \b_0$ and $\n \b^2 \ra \8$}.
\end{equation}
\end{theorem}
This estimate implies (\ref{eq:decay*of*b_c2}) 
as well as (\ref{n_c<8}):
\Corollary{decay*of*b_c}
In the notations of \Thm{th:psi}, we have: 
\\
(i) For all dimension $d$,  $\nu_c <\8$.\\
(ii) For $d \geq 3$, $|\b^\pm_c (d,\n)| \asymp 1/\sqrt{\n}$ as $\n \nearrow \8$.
\end{corollary}
Proof: 
(i) It is enough to show that $p(\beta, \nu )<\nu\lambda$ for negative $\beta$ and large $\nu$.
For fix $\beta <0$ and  $\nu \to \8$, we have 
$ p(\b,\n) \sim \beta \nu$ by \Thm{boundonp}.  In addition to $\b <  \lambda(\b)$, 
this shows that $p(\beta, \nu )<\nu\lambda$  for large $\nu$, and then $\n_c$ is finite.

(ii) Since 
$$
|\b^\pm_c (d,\n)| \st{\mbox{\scriptsize 
(\ref{bc*sand*rough})}}{\ge} |\ln (1\pm \sqrt{a_{L^2}/\n})|
\sim \sqrt{a_{L^2}/\n},
$$ 
we concentrate on the upper bound. 
We first consider the limit $\b \ra 0$ and $\n \b^2 \ra \8$. 
Note that $\n \lm =\n \b +\n \b^2/2 +o (\n \b^2)$ in this limit. 
Comparing this to (\ref{eq:calO(...)}), 
we see that there exist 
a small $\b_0 >0$ and a large $M>0$ such that
$p(\nu, \beta ) < \n \lm$ if $|\b| \le \b_0$ and $\n \b^2 \ge M$. 
Hence, by monotonicity (\Lem{mono*nu}(b)--(c)), we have 
$p(\nu, \beta ) < \n \lm$, if $\n \ge (M /\b_0^2) \vee (M /\b^2)$. 
This implies $\n \b^\pm_c (d,\n)^2 \le M$ if $\n \ge M /\b_0^2$, 
which finishes the proof. 
\hfill $\Box$

\vvs
With \Thm{boundonp} at hand, we can complete: 

\vvs
\noindent {Proof of Theorem \ref{Thm.th:pathlocfull}}: 
By \Lem{lem:2to1}, it is enough to prove (\ref{eq:pathlocfull1}). 
We assume $\b>0$, the other case being similar. Recall (\ref{psi=lm-vp}).
By the convexity of $\widehat{p}_t (\b )$ in $\b$, 
$$
\widehat{p}_t (2\b )-\widehat{p}_t (\b ) 
 \ge  \b {\6 \widehat{p}_t (\b ) \over \6 \b} 
\st{\mbox{\scriptsize (\ref{3100})}}{=}
 {\b \n  \lm (\b) \over t}\int_{[0,t]\times \rd} dsdx
\; Q \; {\mu_t( \chi_{s,x})-\mu_t( \chi_{s,x})^2  
\over 1+ \lm \mu_t( \chi_{s,x})}.
$$
Note that ${1+\lm (\b) u} \le e^{\b} \le e^{\b_0}$ if $0 \le u \le 1$.
Therefore,    
\bdnn
e^{\b_0}\;{\widehat{p}_t (2\b )-\widehat{p}_t (\b ) \over \b \n  \lm (\b)} 
& \ge  & {1\over t}\int_{[0,t]\times \rd} dsdx 
Q  \lef[ \mu_t( \chi_{s,x}) -\mu_t( \chi_{s,x})^2\ri] \\
& = & 1-Q\mu_t^{\otimes 2} (R_t),
\ednn
and hence,
$$
1-\inflim_{t \ra \8}Q\mu_t^{\otimes 2} (R_t)
\le e^{\b_0}{\widehat{p} (2\b )-\widehat{p} (\b ) \over \n  \lm \b} 
\st{\mbox{\scriptsize (\ref{eq:calO(...)})}}{=}
\cO \lef( (\n \b^2)^{-1/6}\ri).
$$
This is (\ref{eq:pathlocfull1}). Now, claim (\ref{eq:pathlocfull2}) 
is (\ref{eq:2to1-2}) in the two-to-one lemma \Lem{lem:2to1}.
\hfill $\Box$

\vvs
\noindent {Proof of Theorem \ref{Thm.th:pathlocfull2}}: 
We suppose that $|\b| \le \b_0$ and $\n \b^2 \ra \8$. 
It follows from (\ref{eq:pathlocfull1}) that
$$
\limsup_{t\to\8}Q \mu_t^{\otimes 2} (1-R_t)=\cO \lef( (\nu \b^2)^{-1/6}\ri).
$$
Thus,  (\ref{eq:r2to1-3})--(\ref{eq:r2to1-5}) directly follow 
from this via (\ref{eq:2to1-3})--(\ref{eq:2to1-5}).
\hfill $\Box$
\SSC{Bound on $p(\b, \n)$,  proof of \Thm{boundonp}}
This section is devoted to the proof of the bound in \Thm{boundonp}.
Before going into the technical details, we sketch the 
\subsection{General strategy}
We first explain the strategy of the estimate in the 
regime $|\b| \le \b_0 <\8$ and $\n \b^2 \ra \8$. 
Since $\nu \b$ diverges (to $\pm \8$),
it is natural to normalize the partition function and define
\bdnl{Zhat}
\hat Z_t= Z_t \exp\{-\nu \b t\} = P[\exp \lef( \b [\h (V_t) -\nu t]\rig)]
\edn
with a centering of the hamiltonian. 
Now, the aim is to bound  $\hat Z_t$ from above, since 
$p(\nu, \b) \geq \nu \b$ by (\ref{psile}).
With two parameters $\gm, \del >0$ to be fixed later on, 
we will define an event of the environment $\cM_{t,\gm, \del} \sub \cM$  
such that 
\bdnl{QMgm->1}
\lim_{t \nearrow \8}\mbox{$Q[ \cM_{t,\gm, \del} ]$}=1,
\edn
and 
\bdnl{exp(lm^1/2)}
\suplim_{t \to \8} \frac{1}{t}
\ln  Q \left[\hat Z_t\; ;\; \cM_{t,\gm, \del}\right]=
\cO \big((\nu \b^2)^{5/6}\big).
\edn
Then, in order to conclude (\ref{eq:calO(...)}), we first
observe that
\bdmn
|\; Q \left[ \ln Z_t |\cM_{t,\gm, \del}\right]
-
Q \left[  \ln Z_t \right] \;|
&\leq& 
 Q \left[ |  \ln Z_t 
- Q[ \ln Z_t]| 
\;|\; \cM_{t,\gm, \del}\right] \nn \\ 
&\leq& 
 Q \left[ |  \ln Z_t 
- Q[ \ln Z_t]| \right]/Q[\cM_{t,\gm, \del}] \nn \\
&=&
{\cal O}(\sqrt{t}) \label{t^1/2}
\edmn
by the concentration property (2.31) in \cite{CY05} and (\ref{QMgm->1}). 
Therefore, 
we have for fixed $\b, \nu$,
\bdnn
p(\b, \nu) &=& 
\lim_{t \to \8} \frac{1}{t}Q \left[  \ln Z_t \right] \\
&\stackrel{ (\ref{t^1/2})}{=}&
\lim_{t \to \8} \frac{1}{t}Q \left[  \ln Z_t \;|\; \cM_{t,\gm, \del}\right]\\
&\stackrel{\rm Jensen}{\leq }&
\inflim_{t \to \8} 
\frac{1}{t} \ln Q \left[  Z_t \;|\; \cM_{t,\gm, \del}\right]\\
&\stackrel{ (\ref{QMgm->1})}{=}&
\inflim_{t \to \8} \frac{1}{t}
\ln  Q \left[\hat Z_t\; ;\; \cM_{t,\gm, \del}\right] + \nu \b
\ednn
This, together with (\ref{exp(lm^1/2)}) and the first bound in (\ref{psile}), 
proves (\ref{eq:calO(...)}).
\medskip

\subsection{The set of good environments}
Let us now turn to the construction of the set $\cM_{t,\gm, \del}$. 
For $a,t>0$ and continuous paths $\{f, g\} \sub \W$, we denote by 
$\rho_{t}(f,g)=\sup_{s \in [0,t]}\|f(s)-g(s)\|_\8$ the uniform distance on $[0,t]$ for the supremum norm in $\R^d$,
and by $\cK_{t,a}$ the set of absolutely continuous function 
$f: \R \to \R^d$, such that $f(0)=0$ and
$$
\frac{1}{t}\int_0^t |{\dot{f}}(s)|^2 ds \leq a^2\;.
$$ 
\Lemma{cameron-martin}
Let $a,\del \in (0,\8)$ and $a\del \ge 1/2$. 
Then, there exists $t_0(\delta)<\8$ such that for all $t \geq t_0(\delta),$
$$P[ \rho_{t}(B, \cK_{t,a}) \ge \del ] 
\leq 2d \exp \lef( -\mbox{$\frac{a^2t}{2}$}+b(t, a,\del ) \ri), 
$$
where
\begin{equation} \label{def:b}
b(t,a,\del )=2t\del^{-2} \ln (2a\del ) [1+\ln (2a\del )]
+\ln (2a\del) \in (0,\8).
\end{equation}
\end{lemma}
The lemma follows from a result of Goodman and Kuelbs \cite[Lemma 2]{GK91}  (taking there $\lm=a^2t$ and 
 $\e=\del /\sqrt{t}$) in the case when $t=1$, which can also be applied to cover general $t>0$ using 
 the scaling property of Brownian motion.

Now, we construct the event $\cM_{t,\gm, \del}$. 
We first cover 
the compact set $\cK_{t,a}$ with finitely many $\rh_t$-balls with radius 1.
The point here is that the number of the balls we need 
is bounded from above explicitly in terms of $a, \del$ and $t$ as we 
explain now. We will use a result of Birman and 
Solomjak, Theorem 5.2 in \cite{BiSo67} (taking there $p=2, \alpha=1, q=\8, m=1, \om=
1$) which yields a precise estimate of the $\epsilon$-entropy of the unit sphere of Sobolev
spaces for $L^q$-norms: For all $ \del >0$,
the set  $\cK_{1,1}$ can be covered by a number
smaller than $\exp \{C_1/\del\}$ of
$\rho_1$-balls  with radius $\del$, where $C_1=C_1(d) \in (0,\8)$
Since, for $a,t>0$, a map $f \mapsto g$, 
$g(u)=(ta)^{-1}f(ut)$ defines a bijection 
from $\cK_{t,a}$ to $\cK_{1,1}$, it follows that,
we can find $f_i \in \cK_{t,a}, 1 \le i \le  
i_0 \le \exp \{C_1 ta/\del\}$, such that
$$\cK_{t,a} \subset \bigcup_{1 \le i  \le i_0}
\Big\{ f \in \W: \rho_t(f,f_i) \leq \del \Big\}\;.
$$

We define the set $\cM_{t,\gm, \del}$ 
of ``good environments'' by
$$
\cM_{t,\gm, \del} = \bigcap_{i=1}^{i_0}\Big\{ \h \in \cM \; ; \; 
\mbox{$\h ( V_t (f_i)) \le
(1+\gm) \n   t $}\Big\}\;.
$$

 Since $\eta( V_t(f_i))$ has Poisson 
distribution with mean $\n t$, 
we have, by union bound and  Cram{\'e}r's upper bound \cite{DeZe98},
for all $t>0$ and $\gm >0$,
\bdmn \nn
Q[ \cM_{t,\gm, \del}^c ] &\leq & \sum_{i  \leq i_0}
\mbox{$ Q[  \eta( V_t(f_i)) > (1+\gm) \n t ]$}\\
\label{poiss} 
&\leq &
\mbox{$ 
\exp\{ -t[\n \lm^*(1+\gamma) - C_1 a/\del]\}$}\;,
\edmn
where 
$$\lm^*(u)=\sup_{\b \in \R} \{ \b u -\lm (\b)\}
=u \ln u -u+1, \qquad u>0.$$
We will eventually take a vanishing $\gamma >0$, so that $\lm^*(1+\gamma) \sim \gamma^2/2$.

\subsection{The moment estimate (\ref{exp(lm^1/2)})}
Letting 
$$\hat \zeta_t=\hat \zeta_t(B)= \exp \{\b \h(V_t)-\n \b  t\},$$ 
we decompose the expectation of $\hat Z_t=P[ \hat \zeta_t]$
on the set $\cM_{t,\gm, \del}$ into
\bdmn \nn
Q \Big[ \hat Z_t ; \cM_{t,\gm, \del}\Big] 
&=& 
Q \Big[ P[ \hat \zeta_t ; \rh_t(B, \cK_{t,a})\! <\! \del] ; 
\cM_{t,\gm, \del}\Big] +
Q \Big[ P[ \hat \zeta_t ; \rh_t(B, \cK_{t,a})\!\geq\! \del] ; 
\cM_{t,\gm, \del}\Big]
\\ \nn
&\leq &
Q \Big[ P[ \hat \zeta_t ; \rh_t(B, \cK_{t,a})\! <\! \del] ; 
\cM_{t,\gm, \del}\Big] +
P \Big[ Q[ \hat \zeta_t] ; \mbox{$ \rh_t(B, \cK_{t,a})\! \geq \! \del$}  \Big]\\
& \leq & P \Big[ Q[ \hat \zeta_t ; \cM_{t,\gm, \del}];  \rh_t(B, \cK_{t,a})\! <\! \del \Big] 
+ 2d \exp \lef( t g (\b,\n, a,\del ) +\ln (2a \del )\ri)\; \label{alali1}
\edmn
using \Lem{cameron-martin} for large enough $t$, and the notation
$$
\hat \lm = \lm - \b \ge 0, 
\; \; \mbox{and}\; \; 
g (\b,\n, \del,a )=\n \hat \lm-\frac{a^2}{2}
+{4 \over \del^2} \ln^2 (2a \del ).
$$
For a path $B$ such that 
$\rh_t(B, \cK_{t,a})\! <\! \del$, which 
contributes to the 
first term in the right-hand side of (\ref{alali1}), 
we can select $i^*=i^*(B) \in \{1,\ldots, i_0\}$
such that $\rh_t(B,f_{i^*})\leq 2 \del$, allowing us to bound the 
$Q$-expectation as
\bdnn 
Q[ \hat \zeta_t(B) ; \cM_{t,\gm, \del}]
&= & Q[ e^{\b[ \h(V_t(B))-\h(V_t(f_{i^*}))]} \times
\hat \zeta_t(f_{i^*}) ; \cM_{t,\gm, \del}] \\ 
&\leq & 
 Q[ e^{\b[ \h(V_t(B))-\h(V_t(f_{i^*}))]}; \cM_{t,\gm, \del}]
 \times e^{\b^+ \n \gamma t}\\ 
&\leq & Q[ e^{\b[ \h(V_t(B))-\h(V_t(f_{i^*}))]}]\times e^{\b^+ \n \gamma t}
\ednn
by definition of  the set  $ \cM_{t,\gm, \del}$ of good environments. 
Now observe that the exponent only collects a few points from the Poissonian environment
(this is where we significantly improve on the estimate in \cite{CY05}):
For paths $B$ such that $\rh_t(B, \cK_{t,a})\! <\! \del$, 
$$
\h(V_t(B))-\h(V_t(f_{i^*})) = 
\h\Big( V_t(B) \setminus V_t(f_{i^*})\Big) - \h 
\Big( V_t(f_{i^*}) \setminus V_t(B)\Big),
$$
that is, the difference of two independent Poisson variables
with the same parameter given by 
$u=\nu |V_t(B) \setminus V_t(f_{i^*})| $ 
(with $|\cdot|$ denoting here the Lebesgue measure in $\R^{1+d}$). 
By simple geometric considerations, we see that
this parameter is bounded by $C_2 t \nu \del $ with 
a constant $C_2=C_2(d)$, 
and the above expectation can be bounded using this parameter.  
Precisely, we obtain for such $B$'s,
\bdnn
Q[ e^{\b[ \h(V_t(B))-\h(V_t(f_{i^*}))]}] &=&
\exp\{ u(\lm(\b)+\lm(-\b))\}\\ 
&\le &
\exp\{ C_2t\nu\del (\lm(\b)+\lm(-\b))\} \\
&\le & \exp ( C_3t\nu \del \b^2),
\end{eqnarray*}
where $C_3=C_3(d,\b_0)$. 
Hence, 
\bdnl{eq:finnn}
 P \Big[ Q[ \hat \zeta_t(B) ; \cM_{t,\gm, \del}];  
\rh_t(B, \cK_{t,a})\! <\! \del \Big] 
\le
\exp ( th (\b,\n,\del,\gm ) ), 
\edn 
where 
$h (\b,\n,\del,\gm )=C_3\nu \del \b^2+\b^+\n\gamma$.
Therefore, by (\ref{alali1}) and (\ref{eq:finnn}),
\bdnl{exp(lm^1/2)<}
\suplim_{t \to \8} \frac{1}{t}
\ln  Q \left[\hat Z_t\; ;\; \cM_{t,\gm, \del}\right]
\le g (\b,\n, \del,a ) \vee h (\b,\n,\del,\gm ).
\edn
To tune the parameters, we set\footnote{We write 
$b(\b,\n) \asymp b'(\b,\n)$ if there exist positive 
finite constants $C_-, C_+$ such that $C_-b \leq b' \leq C_+b$ for 
$|\b| \leq \b_0, \n \geq 1$.} 
\bdnn
a & = & ( 2\n \hat \lm)^{1/2} \asymp (\n \b^2)^{1/2} \ra \8, \\
\del & = & a^{-1/3}\asymp (\n \b^2)^{-1/6}\ra 0, \\
\gamma & = & \lef( {4C_1 a \over \n \del} \ri)^{1/2} 
\asymp |\b|(\n \b^2)^{-1/6} \ra 0.
\ednn
Note that $a \del =a^{2/3} \ra \8$ (hence \Lem{cameron-martin} is available),  
and for $\n \b^2$ large enough,
$$
\n \lm^*(1+\gm)\geq {\n \gm^2 \over 3}={4C_1 a \over 3 \del},
$$
so that (\ref{QMgm->1}) is satisfied by (\ref{poiss}). Moreover,
$$
g (\b,\n, \del, a) 
 \asymp (\n \b^2)^{1/3}\ln^2 (\n \b^2), \; \; \; 
h (\b,\n,\del,\gm )  \asymp  (\n \b^2)^{5/6}.
$$
Thus, we get (\ref{exp(lm^1/2)}) from (\ref{exp(lm^1/2)<}).
This ends the proof.
\hfill $\Box$

\medskip
\SSC{Estimates of the critical curves}
In this section we prove the results for the phase diagram.
\subsection{Some auxiliary curves and the proof of \Thm{bc}} 
\label{auxiliary}
We first remark that the function $p(\b,\n)$ is
monotone and smooth in both variables $\b, \n$. 
\Lemma{mono*nu}
\bds
\item[(a)]
\bdnl{mono*nu1}
\b (\n-\n') \leq p(\b,\n) -p(\b,\n')\leq   \lm(\b)(\n-\n')
\; \; \; \mbox{for $0<\n'<\n$ and $\b \in \R$.}
\edn
\item[(b)]
$(0,\8)\times [0,\8) \ni (\nu, \b) \mapsto \lm (\b)\nu -p(\b,\n)$ 
is non-decreasing in both $\nu$ and $\b$.
\item[(c)]
$(0,\8)\times (-\8,0] \ni (\nu, \b) \mapsto \lm (\b)\nu -p(\b,\n)$ 
is non-decreasing in $\nu$ and non-increasing in $\b$.
\item[(d)]
$(0,\8)\times \R \ni (\nu, \b) \mapsto p(\b,\n)$ is continuous. 
\eds 
\end{lemma}
Proof: 
(a) 
On an enlarged probability space, 
we can couple $\h$ with two mutually independent Poisson point 
processes $\h', 
\h''$ with intensity measures $\n' dsdx, 
(\n-\n')dsdx$, so that $\h = \h' + \h''$. 
Denote by $Q', Q''$ the expectation with respect to 
these new Poisson processes.
By Jensen inequality and by independence, we have,
$$
Q \ln Z_t=Q'Q''\ln P[ \exp \{ \b \h (V_t) \}]
\lef\{ \barray{l}
\ge Q' \ln P[ \exp \{ \b Q'' \h (V_t) \}], \\
\le Q' \ln Q''P[ \exp \{ \b  \h (V_t) \}] ,
\earray \ri.
$$
leading to 
$$
p_t(\b,\n') +\b (\n-\n') \leq p_t(\b,\n) 
\leq p_t(\b,\n') + \lm(\b)(\n-\n'), 
$$
cf. (\ref{pt(b,nu)}).
This proves (\ref{mono*nu1}) via (\ref{psias}). \\
(b)--(c):
These follow from (\ref{mono*nu1}) and \Thm{th:psi}(b). \\
(d) We also see from (\ref{mono*nu1}) that 
$\nu \mapsto p(\n, \b)$ is uniformly continuous, when $\b$ varies 
over a compactum. This, together with the continuity of 
$\b \mapsto p(\n, \b)$, shows (d).
\hfill $\Box$

Let us state now the main technical result, which will be proved in the next section.
It reveals that the families of curves
\bdnl{eq:courbeC}
\cC^\a_a:\; \; \; \n(\b)=a |\lm(\b)|^{-\a},\; \; \b \in \R \bsh \{ 0 \} 
\edn
with $\alpha \in [1,\8)$ and $ a >0,$
convey relevant  information on the critical curve. 
\Theorem{th:phasediag}
Assume $d \geq 3$. 
\bds
\item[(a)] 
Let $0<\b_0<\b_1$ with $(\b_0,\n_0) \in \Crit$ 
(Hence $\b_0=\b_c^+(\n )$). Then, for $\b>0$, 
\bdmn
\n > \n_0,\; \n \lm(\b)^2 > \n_0 \lm(\b_0)^2 & \Rightarrow &
(\b,\n) \in \cL , \label{crit*sand+1}\\ 
0 < \b \le \b_0, \; \; 0< \a \le \a (\b_0), \; \; 
\n \lm(\b)^\a  \leq \n_0 \lm(\b_0)^\a
 &\Rightarrow& (\b,\n)\in \cD , \label{crit*sand+2}\\ 
\n \le \n_0,\; \;\n \lm(\b)^2 \le \n_0 \lm(\b_0)^2 &\Rightarrow& 
(\b,\n) \in \cD,
\label{crit*sand+3}\\ 
\b_0 < \b \le \b_1, \; \; 0< \a \le \a (\b_1), 
\; \; \n \lm(\b)^\a > \n_0 \lm(\b_0)^\a
 &\Rightarrow&
(\b,\n) \in \cL.
\label{crit*sand+4}
\end{eqnarray}
\item[(b)] 
Let $\b_1<\b_0<0$ with $(\b_0,\n_0) \in \Crit$ 
(Hence $\b_0=\b_c^-(\n )$). Then, for $\b <0$, 
\bdmn
\n >\n_0,\; \; \a (\b_0) \le \a, \; \; 
 \n |\lm(\b)|^{\a}  > \n_0 |\lm(\b_0)|^{\a} &\Rightarrow& 
 (\b,\n)\in \cL, 
\label{crit*sand-1}\\
 \b_0 \le \b, \; \;
\n \lm(\b)^2 \le \n_0 \lm(\b_0)^2
&\Rightarrow& (\b,\n) \in \cD , \label{crit*sand-2}\\ 
\b_1 \le  \b < \b_0, \; \; \a (\b_1)  \le \a, \; \; 
\n |\lm(\b)|^\a  \leq \n_0 |\lm(\b_0)|^\a
 &\Rightarrow&
(\b,\n)\in \cD, \label{crit*sand-3}\\
\b < \b_0, \; \; \n \lm(\b)^2 > \n_0 \lm(\b_0)^2
&\Rightarrow& 
(\b,\n) \in  \cL.
\label{crit*sand-4}
\end{eqnarray}
\eds
\end{theorem}
Note that, since $\alpha$ is monotone, it would have been sufficient to state the above 
results taking 
$\alpha=\alpha(\beta_0)$ in (\ref{crit*sand+2}), (\ref{crit*sand-3}), and 
taking 
$\alpha=\alpha(\beta_1)$ in (\ref{crit*sand+4}), (\ref{crit*sand-3}).

We summarize the results in figure \ref{f-diag}.
%
%
\begin{figure} \centering
\caption{Estimates on the critical curve, $d \geq 3$.} \label{f-diag2}
\end{figure}

\medskip

\vvs 
\noindent Proof of \Thm{bc}.
The upper and the lower bounds of $\b^+_c (\n)$ 
in (\ref{bc+sand+1}) follow from 
(\ref{crit*sand+1}) and (\ref{crit*sand+2}), respectively. 
Similarly, 
(\ref{crit*sand+3})--(\ref{crit*sand+4}) imply (\ref{bc+sand-1}), 
(\ref{crit*sand-1})--(\ref{crit*sand-2}) imply (\ref{bc-sand+1}), 
(\ref{crit*sand-3})--(\ref{crit*sand-4}) imply (\ref{bc-sand-1}).
Noting that 
$$\b^\pm_c (\n_0)=\ln (1\pm c^\pm_1 (\n_0))= - \ln (1 \mp c^\pm_2 (\n_0)),
\qquad 
 c^\pm_2 (\n_0)=\frac{c^\pm_1 (\n_0)}{1 \pm c^\pm_1 (\n_0)},
 $$
(\ref{bc+sand+2}),(\ref{bc+sand-2}),(\ref{bc-sand+2}),(\ref{bc-sand-2}) 
follow easily from 
(\ref{bc+sand+1}),(\ref{bc+sand-1}),(\ref{bc-sand+1}),(\ref{bc-sand-1}), 
respectively. We give the proof of  (\ref{bc+sand+2}): by convexity,
\begin{equation} \label{eq:logconv}
\frac{1}{1+x} \leq \frac{\ln (1+x)-\ln(1+y)}{x-y} \leq \frac{1}{1+y} , \quad -1<y<x.
\end{equation}
By the upper bound in (\ref{bc+sand+1}), 
\begin{eqnarray*}
\b_c^+(\n_0)-\b_c^+(\n) &\geq& 
\ln \lef( 1+c_1^+(\n_0)\ri) -
\ln \lef( 1+c_1^+(\n_0)\lef( \n_0 \over \n\ri)^{1/2}\ri)\\
& \stackrel{(\ref{eq:logconv})}{\geq} & \frac{c_1^+(\n_0)}{1+c_1^+(\n_0)} 
 \lef( 1-\lef( \n_0 \over \n\ri)^{1/2} \ri) ,
\end{eqnarray*}
which is the lower bound in (\ref{bc+sand+2}). To prove the other one, 
we use the lower bound in (\ref{bc+sand+1}), 
\begin{eqnarray*}
\b_c^+(\n_0)-\b_c^+(\n) &\leq&
\ln \lef( 1+c_1^+(\n_0)\ri) -
\ln \lef( 1+c_1^+(\n_0)\lef( \n_0 \over \n\ri)^{1/\alpha}\ri) 
\\
& \stackrel{(\ref{eq:logconv})}{\leq} & 
\frac{c_1^+(\n_0)}{1+c_1^+(\n_0)\lef( \n_0 \over \n\ri)^{1/\alpha} }
 \lef( 1-\lef( \n_0 \over \n\ri)^{1/\alpha} \ri) \\
 &=&
 \frac{c_1^+(\n_0)}{\lef( \n \over \n_0\ri)^{1/\alpha}+c_1^+(\n_0) }
 \lef( \lef( \n \over \n_0\ri)^{1/\alpha} -1\ri) \\
 &\leq & 
\frac{c_1^+(\n_0)}{1+c_1^+(\n_0) }
 \lef( \lef( \n \over \n_0\ri)^{1/\alpha} -1\ri),
\end{eqnarray*}
which is the upper bound in (\ref{bc+sand+2}). The other estimates are left to the reader.
\hfill $\Box$
\subsection{Proof of \Thm{th:phasediag}}
\label{p*th:phasediag}

\vvs
The following lemma plays an important role in the proof of 
\Thm{th:phasediag}. 
\Lemma{dp/dn}
We have 
\begin{equation}
\label{eq:dnp}
t\;\frac{\6 p_t}{\6 \n} 
(\beta,\nu)  =\int_{[0,t]\times \rd} dsdx
\; Q \;  \ln [1+\lm \mu_t( \chi_{s,x})] \;.
\end{equation}
\end{lemma}
Proof:
For $k \geq 1$, let 
$$
Z_{t,k}=P[ e^{\b \h(V_t)}; A_k]\;,
\; \; \; 
\mbox{with} \; \; A_k=\{B_s \in [-k,k]^d, \;\forall s \leq t\},
$$
and $p_{t,k}(\b,\n )=t^{-1} Q \ln Z_{t,k}$. 
With $r=r_d$ the radius of $U(0)$, we let $K^r=[-k-r,k+r]^d$.
In this proof, we write $Q=Q_\n$, and we use that, 
when restricted to the bounded set $K_t=(0,t]\times K^r$,  
the Poisson point measure $Q_\n$ is absolutely continuous with respect to 
the Poisson point measure $Q_1$ with unit intensity, in order to write
$$
Q_\n \ln Z_{t,k}= Q_1 \lef[ \rh_{t,\nu} \ln Z_{t,k}\ri] 
\; \; \;\mbox{with}\; \; \; 
\rh_{t,\nu}={d Q_\n |_{K_t}\over d Q_1 |_{K_t}}=
\exp \lef( \h_t(K_t)\ln \n 
-  (\n-1)t|K^r|\ri).
$$
Thus, $t p_{t,k}(\b,\n)$ is differentiable in $\nu$, with derivative
\begin{eqnarray*}
Q_1 \lef[ \rh_{t,\nu}
\; \frac{\h_t(K_t) \!-\!  \n t|K^r|}{\n}  \ln Z_{t,k} \ri]
&=& 
\frac{1}{\n} Q_\n 
\lef[  \int_{K_t}
 \widehat{\h}_t(dsdx)   \ln Z_{t,k}\ri]\\
&\st{\mbox{\scriptsize (\ref{ipp1})}}{=}& 
\int_{K_t} dsdx \; Q_\n     \ln 
\frac{Z_{t,k}(\h_t+\del_{s,x})}{Z_{t,k}(\h_t) }  \\
&=& 
\int_{[0,t]\times \rd} dsdx \; Q_\n   \ln 
\frac{Z_{t,k}(\h_t+\del_{s,x})}{Z_{t,k}(\h_t) }\\
&=& 
\int_{[0,t]\times \rd} dsdx \; Q_\n   
\ln \lef( 1+\lm \m_t (\chi_{s,x} | A_k)\ri),
\end{eqnarray*}
where the last equality is obtained similarly 
as we did in (\ref{6p/6b}). Now, we write
\bdnl{ptk-ptk}
tp_{t,k}(\b,\n)-tp_{t,k}(\b,1)
= \int_1^\n d\n' \int_{[0,t]\times \rd} dsdx \; Q_{\n'}   \ln 
\lef( 1+\lm \m_t (\chi_{s,x} | A_k)\ri).
\edn
We infer from this identity that
\bdnl{pt-pt}
t p_t(\b,\n)-t p_t(\b,1)= \int_1^\n d\n' 
\int_{[0,t]\times \rd} dsdx \; Q_{\n'}   \ln 
\lef( 1+\lm \mu_t(\chi_{s,x})\ri) \;,
\edn
which shows that $p_t(\b, \cdot)$ 
is differentiable, and also yields (\ref{eq:dnp}).
It is clear that 
$$
p_{t,k}(\b,\n) 
\st{k \nearrow \8}{\to} p_t(\b,\n).
$$
To show that the right-hand side of (\ref{ptk-ptk}) converges to that 
of (\ref{pt-pt}), it is enough to verify that
\bdnl{ptk->pt}
\lim_{k \ra \8}\int_{\rd}  \; Q_{\n'}   \ln 
\lef( 1+\lm \m_t (\chi_{s,x} | A_k)\ri)dx
=\int_{\rd}  \; Q_{\n'}   \ln 
\lef( 1+\lm \m_t (\chi_{s,x})\ri)dx
\edn
and that
\bdnl{ptk<}
\lm -{|\lm|^2 \over 2}
\le
\int_{\rd}  \; Q_{\n'}   \ln 
\lef( 1+\lm \m_t (\chi_{s,x} | A_k)\ri)dx \le \lm.
\edn
After these, it only remains to apply the bounded convergence theorem to 
take care of the rest of the integrations.
The bound (\ref{ptk<}) follows from the elementary inequality:
$$
u-{u^2 \over 2} \le \ln (1+u) \le u, \; \; \; u>-1.
$$
On the other hand, we have that 
\bdnn
\lefteqn{e^{-|\b|}
|\ln (1+\lm \m_t (\chi_{s,x} | A_k))-\ln (1+\lm \m_t (\chi_{s,x}))| } \\
& \st{\mbox{\scriptsize (\ref{eq:logconv})}}{\le} & 
|\m_t (\chi_{s,x} | A_k)-\m_t (\chi_{s,x})| \\
&\le & \m_t (A_k)^{-1}\m_t (\chi_{s,x} ; A_k^{\complement})
+\m_t (A_k)^{-1}\m_t (\chi_{s,x})\m_t (A_k^{\complement}).
\ednn
and hence that
$$
e^{-|\b|}\int_{\rd}\lef| \ln \lef( 1+\lm \m_t (\chi_{s,x} | A_k)\ri)
-\ln (1+\m_t (\chi_{s,x})) \ri| dx
\le 2\m_t (A_k)^{-1}\m_t (A_k^{\complement})
 \st{k \ra \8}{\lra}  0,
$$ 
which proves (\ref{ptk->pt}).
\hfill $\Box$

\vvs 
We will also need the following elementary observation.
\Lemma{elem}
Let
\bdnl{h*alpha}
h_\a (u)=\ln (1+u) -u+{u^2 \over \a (1+u)}, 
\; \; \; \a>0, \; \; u >-1.
\edn 
Then,
\bdnl{ele-cal}
h_\a (u)\lef\{ \barray{ll}
\le 0 & \mbox{if $\a =2$ and $u \in [0,\8)$}, \\
\ge 0 & \mbox{if $\a =2$ and $u \in (-1,0]$}, \\
\ge 0 & \mbox{if $\b>0$, $\a \le \a (\b)$ and $u \in [0,\lm ]$}, \\
\le 0 & \mbox{if $\b<0$, $\a (\b) \le \a $ and $u \in [\lm ,0]$}.
\earray \rig.
\edn 
\end{lemma}
Proof:
We have $h_\a' (u)= {u \lef( 2-\a -(\a-1)u\ri) \over \a (1+u)^2}$.
Therefore, if $1<\a \le 2$, 
\bdnl{ha+}
h_\a \lef\{ \barray{l} 
\mbox{decreases from $\8$ to 0 on $(-1,0]$}, \\
\mbox{increases from $0$ on $I_0\st{\rm def}{=}(0,{2 -\a \over \a -1}]$, 
($I_0=\epty$, when $\a=2$)}, \\
\mbox{decreases on $I_+\st{\rm def}{=}[{2 -\a \over \a -1},\8)$}. \earray \ri.
\edn
and if $\a \ge 2$,
\bdnl{ha-}
h_\a \lef\{ \barray{l} 
\mbox{decreases from $\8$ 
on $J_-\st{\rm def}{=}(-1,{2 -\a \over \a -1}]$}, \\
\mbox{increases to 0 on $J_0\st{\rm def}{=}({2 -\a \over \a -1},0]$, 
($J_0=\epty$, when $\a=2$)}, \\
\mbox{decreases on $[0, \8)$}. \earray \ri.
\edn
The first two lines of (\ref{ele-cal}) follow immediately from 
either (\ref{ha+}) or (\ref{ha-}). 
By the obvious monotonicity, we may assume that 
$\a>1$ to prove the third line of (\ref{ele-cal}). 
Suppose that $\b>0$. Then, 
$\a (\b) < 2$ and hence, for $\a>1$, 
$$
\a \le \a (\b)
\; \; \st{\mbox{\scriptsize (\ref{alpha(b)})}}{\Llra} \; \; h_{\a}(\lm) \ge 0
\; \; \st{\mbox{\scriptsize (\ref{ha+})}}{\Llra} \; \; \min_{[0,\lm]}h_{\a} \ge 0.
$$
This proves the third line of (\ref{ele-cal}) for $\a>1$, 
which proves the case of $0<\a \le 1$ by the obvious 
monotonicity. 
Suppose on the other hand that $\b<0$. Then, 
$\a (\b) > 2$ and hence,
$$
\a \ge \a (\b)
\; \; \st{\mbox{\scriptsize (\ref{alpha(b)})}}{\Llra} \; \; h_{\a}(\lm) \le 0
\; \; \st{\mbox{\scriptsize (\ref{ha-})}}{\Llra} \; \; \max_{[\lm,0]}h_{\a} \le 0.
$$
This proves the fourth line of (\ref{ele-cal})
\hfill $\Box$

\vvs
\noindent Proof of \Thm{th:phasediag}. 
(a) We look for a smooth function 
$\n =\n (\b)$, ($\b \in \R \bsh \{ 0 \})$ such that 
\bdnl{C:dec/inc}
\mbox{$\n (\cdot)$ is decreasing on $(0,\8)$, and is increasing on $(-\8,0)$}.
\edn 
and 
$$
F(\b)\st{\rm def.}{=} \lef.  (tp_t(\b,\n)-t\n\lm(\b) )\ri|_{\nu =\nu (\b)}
$$
is monotone.
We set $G(\b,\n)=tp_t(\b,\n)-t\n\lm(\b)$ and we compute 
with (\ref{eq:dnp}) and (\ref{3101}), 
\bdmn
F'(\b) &=& 
\n' \frac{\6 G}{\6 \n}  (\b,\nu(\b))
+   \frac{\6 G}{\6 \b}  (\b,\nu(\b)) \nn \\
&=& \n'  Q \int_{[0,t]\times \rd} dsdx
\left\{  \ln [1\!+\!\lm \mu_t( \chi_{s,x})]\! 
-\! \lm \mu_t( \chi_{s,x}) \right\} \nn \\
& &  - \n e^\b \lm 
 Q \int_{[0,t]\times \rd} dsdx 
\frac{\mu_t( \chi_{s,x})^2}{ 1+\lm \mu_t( \chi_{s,x})} \nn \\ 
&=& \n' Q \int_{[0,t]\times \rd} dsdx
\left\{  \ln (1\!+\!\lm \mu_t( \chi_{s,x}))\! -\! \lm \mu_t( \chi_{s,x})  
- \frac{\n}{\n'} e^\b \lm 
 \frac{\mu_t( \chi_{s,x})^2}{ 1+\lm \mu_t( \chi_{s,x})}\right\},
\label{eq:tired}
\edmn
if $\n'(\beta)$ does not vanish. We now take the curve 
(\ref{eq:courbeC}) for which (\ref{C:dec/inc}) is satisfied.
We then obtain
\bdnl{F'}
F'(\b) = \n' Q \int_{[0,t]\times \rd} h_\a (\lm \mu_t( \chi_{s,x}))dsdx, 
\; \; \; \mbox{cf. (\ref{h*alpha})}.
\edn
We first prove (\ref{crit*sand+2}). 
By \Lem{mono*nu}(b), it is enough to show that
\bdnl{crit*sand+2=}
0 < \b \le \b_0, \; \; 0< \a \le \a (\b_0), \; \; 
\n \lm(\b)^\a  = \n_0 \lm(\b_0)^\a 
\Rightarrow (\b,\n)\in \cD.
\edn
We take $\a \le \a (\b_0)$ (hence $\a \le \a (\b)$) 
and $a=\n_0 \lm(\b_0)^\a$ in (\ref{eq:courbeC}). 
We then use (\ref{ele-cal}) (the third line), 
(\ref{F'}) and (\ref{C:dec/inc}) to see that 
$F|_{[0,\8)}$ is non-increasing, and hence that 
$0 \ge F(\b) \ge F(\b_0)=0$ if $\b \in [0, \b_0]$. 
This proves (\ref{crit*sand+2=}).\\
The proof of (\ref{crit*sand+3}) is similar as above. 
Since $(0,\b_0] \times (0,\n_0] \sub \cD$ by the monotonicity 
(\Lem{mono*nu}(b)), we may assume $\b_0 \le \b$. 
By \Lem{mono*nu}(b) again, it is enough to show that
\bdnl{crit*sand+3=}
\b_0 \le \b,\;\;\n \lm(\b)^2 = \n_0 \lm(\b_0)^2 \Rightarrow 
(\b,\n) \in \cD.
\edn 
We take $\a =2$ and $a=\n_0 \lm(\b_0)^2$ in (\ref{eq:courbeC}).
We then see from 
(\ref{ele-cal}) (the first line), (\ref{F'}) and (\ref{C:dec/inc}) that 
$F|_{[0,\8)}$ is non-decreasing, and hence that 
$0 = F(\b_0) \le F(\b) \le 0$ if $\b_0 \le \b$. 
This proves (\ref{crit*sand+3=}). \\
We next prove (\ref{crit*sand+1}). 
Choose $\e>0$ such that $\n \ge \n_0 +\e$ and 
$\n \lm(\b)^2  \ge (\n_0 +\e)\lm (\b_0+\e)^2$. 
Then, $(\n_0+\e, \b_0+\e) \in \cL$ 
by \Lem{mono*nu}(b) and the fact that $(\n_0, \b_0) \in \6 \cL$. 
Since $[\b_0+\e,\8 ) \times [\n_0+\e,\8 ) \sub \cL$ by the monotonicity 
(\Lem{mono*nu}(b)), we may assume $\b \le \b_0+\e$. 
By \Lem{mono*nu}(b) again, it is enough to 
show that
\bdnl{crit*sand+1=}
0<\b \le \b_0+\e,\; \n \lm(\b)^2 = (\n_0 +\e)\lm (\b_0+\e)^2  \Rightarrow 
(\b,\n) \in \cL.
\edn
We take $\a =2$ and $a=(\n_0 +\e)\lm (\b_0+\e)^2$ in (\ref{eq:courbeC}).
We then see from 
(\ref{ele-cal}) (the first line), (\ref{F'}) and (\ref{C:dec/inc}) that 
$F|_{[0,\8)}$ is non-decreasing, and hence that 
$F(\b) \le F(\b_0+\e)<0$ for $\b \in (0,\b_0+\e]$. 
This proves (\ref{crit*sand+1=}). \\
The proof of (\ref{crit*sand+4}) is similar as above. 
Choose $\e>0$ such that $\b_0 +\e \le \b$ and 
$\n \lm (\b)^\a \ge (\n_0 +\e)\lm (\b_0+\e)^\a$. 
Then, $(\n_0+\e, \b_0+\e) \in \cL$ 
by \Lem{mono*nu}(b) and the fact that $(\n_0, \b_0) \in \6 \cL$. 
By \Lem{mono*nu}(b) again, it is enough to 
show that
\bdnl{crit*sand+4=}
\b_0 < \b \le \b_1, \; \; 0< \a \le \a (\b_1), 
\; \; \n \lm(\b)^\a = (\n_0 +\e)\lm (\b_0+\e)^\a \Rightarrow
(\b,\n) \in \cL.
\edn
We take $0<\a \le \a (\b_1)$ (hence $\a \le \a (\b)$) and 
$a=(\n_0 +\e)\lm (\b_0+\e)^\a$ in (\ref{eq:courbeC}).
We then see from 
(\ref{ele-cal}) (the third line), (\ref{F'}) and (\ref{C:dec/inc}) that 
$F|_{[0,\8)}$ is non-increasing, and hence that 
$F(\b) \le F(\b_0+\e)<0$ if $\b \in [\b_0+\e, \b_1]$. 
This proves (\ref{crit*sand+4=}). \\
(b): Proofs of (\ref{crit*sand-1})--(\ref{crit*sand-4}) are similar to 
those of (\ref{crit*sand+1})--(\ref{crit*sand+4}), respectively. 
\hfill $\Box$
\medskip

We conclude this paper by proving Remark \ref{rk:order}.
As shown in Remark 2.2.1 in \cite{CY05},  $p$ is differentiable function of  $\b$
at each point of $\Crit$ with  derivative equal to $\n \lambda'(\b)$. By convexity, the sequence of derivatives of 
$p_t(\b,\n)$ converges to the derivative of the limit. By second order expansion in  (\ref{energy-eq}), this implies that 
$$
\lim_{t\to\8} t^{-1} \int_{[0,t]\times \R^d} ds dx \mu_t(\chi_{s,x})^2 =0.
$$
By second order expansion in  (\ref{eq:dnp}), this implies that 
$\frac{\6 p_t}{\6 \n} 
(\beta,\nu)  $ converges to $\lm$ as $t$ diverges. This implies that, at each, $(\b, \n) \in \Crit $  
when  $d \geq 3$, the differential of
 (the smooth function) 
$p_t(\beta,\nu)$ converges as $t \to \8$ to that of the limit $\n \lm(\b)$.
{\footnotesize

}

\end{document}